\documentclass[11pt,reqno]{amsart}
\usepackage{amssymb,amsmath,amsthm}
\usepackage{color}
\usepackage{enumitem}
\usepackage{fullpage}
\usepackage{url}
\usepackage[Symbol]{upgreek}
\usepackage[unicode,pdfusetitle]{hyperref}
\usepackage[mathscr]{euscript}
\usepackage{tikz-cd}
\usepackage{mathtools}
\usepackage{colonequals}
\newcommand*{\logeq}{\ratio\Longleftrightarrow}

\newtheorem{theorem}{Theorem}
\newtheorem*{theorem*}{Theorem}
\newtheorem{definition}[theorem]{Definition}
\newtheorem*{definition*}{Definition}
\newtheorem{question}[theorem]{Question}
\newtheorem*{question*}{Question}

\newtheorem*{conjecture*}{Conjecture}
\newtheorem{convention}[theorem]{Convention}
\newtheorem*{convention*}{Convention}
\newtheorem{corollary}[theorem]{Corollary}
\newtheorem*{corollary*}{Corollary}
\newtheorem*{remark*}{Remark}
\newtheorem{proposition}[theorem]{Proposition}
\newtheorem*{proposition*}{Proposition}
\newtheorem{lemma}[theorem]{Lemma}
\newtheorem*{lemma*}{Lemma}

\newtheorem*{example*}{Example}
\newtheorem{fact}[theorem]{Fact}
\newtheorem*{fact*}{Fact}

\newtheorem*{claim*}{Claim}
\newtheorem{theoremA}{Theorem}

\numberwithin{theorem}{section}
\numberwithin{equation}{section}

\theoremstyle{definition}
\newtheorem{remark}[theorem]{Remark}

\DeclareMathOperator{\h}{h}

\DeclareMathOperator{\dlarge}{lg}
\DeclareMathOperator{\dsmall}{sm}
\DeclareMathOperator{\nl}{nl}
\DeclareMathOperator{\nt}{nt}

\DeclareMathOperator{\rc}{rc}

\DeclareMathOperator{\res}{res}

\DeclareMathOperator{\elem}{elem}
\DeclareMathOperator{\pow}{pow}
\DeclareMathOperator{\st}{st}

\DeclareMathOperator{\Th}{Th}

\newcommand{\N}{\mathbb{N}}

\newcommand{\Q}{\mathbb{Q}}
\newcommand{\R}{\mathbb{R}}
\newcommand{\T}{\mathbb{T}}
\newcommand{\Z}{\mathbb{Z}}

\newcommand{\cL}{\mathcal L}

\newcommand{\cO}{\mathcal O}

\newcommand{\fm}{\mathfrak{m}}

\newcommand{\upl}{\uplambda}
\newcommand{\Upl}{\Uplambda}
\newcommand{\upo}{\upomega}
\newcommand{\Upo}{\Upomega}

\newcommand{\No}{\mathbf{No}}

\renewcommand{\preceq}{\preccurlyeq}
\renewcommand{\succeq}{\succcurlyeq}
\renewcommand{\geq}{\geqslant}
\renewcommand{\leq}{\leqslant}
\renewcommand{\epsilon}{\varepsilon}

\renewcommand{\d}{\operatorname{d}}
\newcommand{\inv}{^{-1}}

\DeclareFontFamily{OMS}{smallo}{}
\DeclareFontShape{OMS}{smallo}{m}{n}{<->s*[.65]cmsy10}{}
\DeclareSymbolFont{smallo@m}{OMS}{smallo}{m}{n}
\DeclareMathSymbol{\smallo}{\mathord}{smallo@m}{79}

\DeclareFontFamily{U}{fsy}{}
\DeclareFontShape{U}{fsy}{m}{n}{<->s*[.9]psyr}{}
\DeclareSymbolFont{der@m}{U}{fsy}{m}{n}
\DeclareMathSymbol{\der}{\mathord}{der@m}{182}
\DeclareMathSymbol{\derdelta}{\mathord}{der@m}{100}

\author{Elliot Kaplan}
\email{elliot.kaplan@math.su.se}
\title{Model completeness for the differential field of transseries with exponentiation}
\address{Matematiska institutionen, Stockholms universitet, Stockholm, Sweden}
\date{\today}
\begin{document}

\begin{abstract}
We show that the expansion of the valued differential field of logarithmic-exponential transseries by its natural exponential function is model complete.
We give an axiomatization of the theory of this expansion that is effective relative to the theory of the real exponential field, and show that any maximal Hardy field satisfies the same axioms.
We use this to prove a transfer theorem between Hardy fields and transseries.
\end{abstract}

\maketitle

\setcounter{tocdepth}{1}
\tableofcontents

\section*{Introduction}
The differential field $\T$ of logarithmic-exponential transseries is a real closed field extension of $\R$, introduced independently by \'{E}calle~\cite{Eca92} in his solution to the Dulac conjecture and by Dahn and G\"{o}ring~\cite{DG87} in their work on Tarski's problem on real exponentiation.
This structure $\T$ serves as a sort of universal domain for studying the asymptotic behavior of non-oscillating solutions to algebraic differential equations over the reals.
In~\cite{ADH17}, Aschenbrenner, van den Dries, and van der Hoeven showed that the elementary theory of $\T$ as an ordered valued differential field is model complete.

There is a natural exponential $\exp$ on the field $\T$ that makes $\T$ an elementary extension of the real exponential field $\R_{\exp}$~\cite{DMM97}.
Let $\T_{\exp}$ be the expansion of the \emph{differential} field $\T$ by $\exp$.
Then $\T_{\exp}$ is a proper expansion of $\T$ in terms of definability, since $\exp$ is not definable in $\T$; see~\cite[Proposition~16.6.7]{ADH17} and~\cite[Corollary~5.9]{Kap26}.
In the introduction of~\cite{ADH17}, it was asked whether the model completeness result for $\T$ can be extended to $\T_{\exp}$.
In this paper, we provide a positive answer to this question.

\begin{theoremA}[Corollary~\ref{cor:maincor}]\label{thm:mcT}
The elementary theory of $\T_{\exp}$ is model complete in the language of ordered valued differential exponential fields.
\end{theoremA}

Note that in~\cite{ADH17}, the notation $\T_{\exp}$ differs from our use here.
There, $\T_{\exp}$ refers to a specific differential subfield of $\T$ that is actually quite wild, model-theoretically (it defines $\Z$; see~\cite{AD05}).

In~\cite{ADH22B}, the authors showed that the theory of $\T$ is completely axiomatized as follows:
\begin{itemize}
\item $\T$ is an $H$-field, as defined by Aschenbrenner and van den Dries~\cite{AD02}; see Section~\ref{sec:asymp}.
\item $\T$ has \emph{small derivation}:\ the derivative of any infinitesimal element is infinitesimal.
\item $\T$ is \emph{Liouville closed}:\ for each $f \in \T$, there are $g,h \in \T$ with $h\neq 0$ such that $g' = h'/h = f$.
\item $\T$ has the \emph{differential intermediate value property}:\ for any differential polynomial $P$ over $\T$ and $f,g \in \T$ with $P(f)<0<P(g)$, there is $y \in \T$ between $f$ and $g$ with $P(y) = 0$.
\end{itemize}
In particular, the theory of $\T$ is decidable~\cite[Corollary~16.6.3]{ADH17}.
Here we show that $\T_{\exp}$ admits a similar axiomatization:

\begin{theoremA}[Corollaries~\ref{cor:DIVP} and~\ref{cor:maincor}]\label{thm:expaxioms}
Let $K$ be an ordered valued differential exponential field.
Then $K$ is elementarily equivalent to $\T_{\exp}$ if and only if it satisfies the following axioms:
\begin{enumerate}[label=\normalfont(\arabic*), align = right]
\item\label{axiom1} $K$ is an $H$-field with small, nontrivial derivation and $K$ satisfies the differential intermediate value property.
\item\label{axiom2} $\exp$ is an ordered group isomorphism $K\to K^>$ and $\exp(f)' = \exp(f)f'$ for $f \in K$.
\item\label{axiom3} The constant field $C$ of $K$, equipped with the restricted exponential $\exp|_{C}$, is elementarily equivalent to the real exponential field $\R_{\exp}$.
\end{enumerate}
\end{theoremA}

Note that the Liouville closedness axiom is no longer needed:\ for $f \in K$, the differential intermediate value property gives $g \in K$ with $g' = f$~\cite[Lemma~2.1]{ADH22B}, and then $(\exp g)'/(\exp g) = f$ by~\ref{axiom2}.
The only part of this axiomatization that is not known to be effective is~\ref{axiom3}.
Thus, the theory of $\T_{\exp}$ is decidable relative to the theory of $\R_{\exp}$; see~\cite{MW96} for more about decidability of this theory.

As with the theory of $\T$, the differential intermediate value property is a bit of an afterthought.
The axioms actually needed to prove model completeness are the more complicated axioms of \emph{newtonianity and $\upo$-freeness}; see Subsection~\ref{subsec:upofree}.
It is also more convenient to work with the logarithm function instead of the exponential function, as one has much more control over the value groups of valued field extensions that are only closed under taking logarithms; see~\cite[Section~4]{DMM94} for evidence supporting these heuristics.
Accordingly, we study \emph{logarithmic $H$-fields}, a class of $H$-fields equipped with a (not necessarily surjective) logarithm.
We call a logarithmic $H$-field with a surjective logarithm an \emph{exponential $H$-field} and we show that every logarithmic $H$-field has an exponential $H$-field extension with the same constant field.
We further show that every exponential $H$-field can be extended to an $\upo$-free newtonian exponential $H$-field with the same constant field.
Every $\upo$-free newtonian exponential $H$-field is automatically Liouville closed; see Remark~\ref{rem:realclosed-log}.

Let $T_{\log}$ extend the theory of logarithmic $H$-fields by axioms asserting that the constant field $C$ equipped with the restricted logarithm $\log|_{C^>}$ is elementarily equivalent to $\R_{\exp}$ and let $T^{\nt}_{\exp}$ extend the theory of $\upo$-free newtonian exponential $H$-fields by the same constant field axioms.
Using Wilkie's celebrated theorem that $\R_{\exp}$ is model complete~\cite{Wil96}, we show the following:

\begin{theoremA}[Theorem~\ref{thm:mainthm}]\label{thm:modelcompanion}
$T^{\nt}_{\exp}$ is model complete and it is the model companion of $T_{\log}$.
\end{theoremA}

Since $\T_{\exp}\models T^{\nt}_{\exp}$, Theorem~\ref{thm:mcT} is an immediate corollary of Theorem~\ref{thm:modelcompanion}.
In fact, any model of $T^{\nt}_{\exp}$ with small derivation is elementarily equivalent to $\T_{\exp}$ by Corollary~\ref{cor:maincor}, so Theorem~\ref{thm:expaxioms} also follows from Theorem~\ref{thm:modelcompanion} and the results in~\cite{ADH22B}.
Perhaps surprisingly, this entails that all models of $T^{\nt}_{\exp}$ are elementarily equivalent to $\R_{\exp}$ as exponential fields; see Corollary~\ref{cor:modelsT}.

To determine the combinatorial complexity of $T^{\nt}_{\exp}$ and to better understand definable sets in $\T_{\exp}$, we will likely need to prove a quantifier elimination result in some reasonable extended language.
At the end of this paper, we give some ideas about what this language might be.
While we do not yet have such a result, our model completeness proof is robust enough to give some information about definable sets:

\begin{theoremA}[Corollaries~\ref{cor:lomin} and~\ref{cor:stablyembedded}]\label{thm:definable}\
\begin{enumerate}[label=\normalfont(\arabic*), align = right]
\item\label{definable1} For every $X \subseteq \T_{\exp}$ definable with parameters and every $b \in \T_{\exp}\cup \{+\infty\}$, there is $a \in \T_{\exp}$ with $a<b$ such that $(a,b)$ is either contained in or disjoint from $X$.
\item\label{definable2} For every $X \subseteq \T_{\exp}^n$ definable with parameters, the intersection $X \cap \R^n$ is definable in the real exponential field $\R_{\exp}$.
\end{enumerate}
\end{theoremA}

In the literature, the first property above is called \emph{local o-minimality} and the second is \emph{stable embeddedness of the constant field}.

The theory $T^{\nt}_{\exp}$ has other models of interest---any maximal Hardy field is a model of $T^{\nt}_{\exp}$ with small derivation by~\cite{Bou76,ADH24}.
We use this to prove the following transfer theorem:

\begin{theoremA}[Theorem~\ref{thm:transfer}]
Let $H\supseteq \R$ be a henselian Hardy field closed under taking logarithms and let $\iota\colon H\to \T_{\exp}$ be an ordered differential field embedding that preserves logarithms.
Then any maximal Hardy field extending $H$ has the same first-order theory as $\T_{\exp}$ in the language of ordered differential logarithmic fields with parameters from $H$ (interpreted in $\T_{\exp}$ via $\iota$).
\end{theoremA}

In particular, a system of exponential-algebraic differential equations over $H$ has a solution in some Hardy field extension of $H$ if and only if its image under $\iota$ has a solution in $\T_{\exp}$.

Conway's field of surreal numbers equipped with the Kruskal--Gonshor exponential and the Berarducci--Mantova derivation is another model of $T^{\nt}_{\exp}$ by~\cite{BM18,ADH19}.
Using this, we answer the first question posed in~\cite{BEK16} (also asked in~\cite{ADH19}).

\begin{theoremA}[Proposition~\ref{prop:Noembedding}, Corollaries~\ref{cor:HardyNo} and~\ref{cor:elemTNo}]
Every model of $T_{\log}$ with archimedean constant field and small derivation admits an ordered differential logarithmic field embedding into $\No$.
In particular, every Hardy field closed under taking logarithms admits such an embedding.
Moreover, the natural embedding $\T_{\exp}\to \No$ is an elementary embedding of ordered differential exponential fields.
\end{theoremA}

The differential field $\T$ has two other expansions worth mentioning:
\begin{itemize}
\item $\T_{\pow}$ expands $\T$ by the map $(f,r) \mapsto f^r\colon \T^>\times \R \to \T^>$.
\item $\T_{\elem}$ expands $\T_{\exp}$ by the sine and cosine functions restricted to the interval $[-1,1]$.
\end{itemize}
The field $\T_{\pow}$ was studied in~\cite{Kap26}, where analogs of the theorems above were established.
A variant of $\T_{\elem}$, where $\exp$ is also restricted to $[-1,1]$, was studied in~\cite[Chapter~8]{Kap21}, and we give a full treatment of $\T_{\elem}$ in Subsection~\ref{subsec:trig} below.
We show in Theorem~\ref{thm:trigmainthm} and Corollary~\ref{cor:elemTderiv} that many of the theorems above hold for $\T_{\elem}$ as well.

In terms of definability, we have 
\[
\T\subsetneq\T_{\pow}\subsetneq\T_{\exp}\subsetneq\T_{\elem},
\]
where $\subsetneq$ here indicates that one is a proper reduct of the other.
See~\cite[Corollary~5.9]{Kap26} for strictness of the first two reducts; the third follows from Corollary~\ref{cor:stablyembedded}, Theorem~\ref{thm:trigmainthm}\ref{trigmainthmv}, and~\cite{Bia97}.
The fact that $\T_{\elem}$ properly expands $\T_{\exp}$ is a bit superficial:\ it really only manifests in the induced structure on the constant field by Theorem~\ref{thm:trigmainthm}\ref{trigmainthmiii}.
Some of the results on $\T_{\pow}$ in~\cite{Kap26} can be deduced from the results here; see Remark~\ref{rem:pow} for a discussion of this.

\subsection*{Strategy and outline}
Our proof of Theorem~\ref{thm:modelcompanion} closely follows the proof in~\cite{ADH17} that the theory of $\upo$-free newtonian Liouville closed $H$-fields is the model companion of the theory of $H$-fields.
We prove a general proposition (Proposition~\ref{prop:logclosedtest2}) that allows us to upgrade many of the extension and embedding results for $H$-fields in~\cite{ADH17} to results about logarithmic $H$-fields.
Nevertheless, some of the results in~\cite{ADH17} do not admit easy generalizations, even with this tool in hand.
In particular, the extension results in~\cite[Sections~13.4 and~16.1]{ADH17} (realizing cuts in the value group and extensions controlled by asymptotic couples) are quite difficult to generalize.
These extensions are handled in Propositions~\ref{prop:zgap} and~\ref{prop:dtranselem} respectively, with all of Section~\ref{sec:vgext} devoted to the proof of Proposition~\ref{prop:dtranselem}.

In proving Proposition~\ref{prop:dtranselem}, we establish some facts about the set $v(a-K)$ for $K$ an $H$-field and $a$ an integral of an element in $K$.
We show that this set does not change when passing from $K$ to its henselization (Proposition~\ref{prop:hensel}), nor does it change when extending $K$ by another integral $b$, assuming $v(a-K)$ and $v(b-K)$ are suitably different (Proposition~\ref{prop:twoints}).
Both results may be of general use (and both hold in the broader context of $\d$-valued fields; see Remark~\ref{rem:dvalued}).

Once Theorem~\ref{thm:modelcompanion} is established, the other main theorems follow relatively easily (building on the work done in~\cite{ADH17,ADH19,ADH22B,ADH24} of course).
Theorem~\ref{thm:definable}\ref{definable2} needs a different argument than the one given in~\cite{ADH17}, which goes through quantifier elimination.
The argument we give is essentially the same as the one we gave for $\T_{\pow}$ in~\cite{Kap26}.

Section~\ref{sec:prelim} is devoted to preliminaries on ordered fields, valued fields, and differential fields.
In Section~\ref{sec:asymp}, we recall background on $H$-fields, and prove Propositions~\ref{prop:hensel} and~\ref{prop:twoints} on approximating integrals.
Logarithmic and exponential $H$-fields are introduced in Section~\ref{sec:logexpfields} and we study logarithmic $H$-field extensions and logarithmic $H$-field embeddings in Section~\ref{sec:exts}.
In Section~\ref{sec:vgext}, we prove Proposition~\ref{prop:dtranselem} on extensions controlled by the asymptotic couple.
We use these results in Section~\ref{sec:mc} to prove our main theorems about $T^{\nt}_{\exp}$.
We also study expansions by restricted trigonometric functions in this section.
In Section~\ref{sec:examples}, we discuss transseries, Hardy fields, and surreal numbers.
Many natural questions remain open, and we discuss these in Section~\ref{sec:open}.

\subsection*{Acknowledgements, some history of this paper, and AI disclosure}
A version of this paper first appeared in April 2020, and I thank Lou van den Dries, Allen Gehret, and Nigel Pynn-Coates for providing helpful feedback on that version.
Unfortunately, that version contained a gap---I had combined Facts~\ref{fact:smallint} and~\ref{fact:largeint} into an erroneous lemma.
This lemma essentially asserted that when extending a henselian $H$-field $K$ by an element $a \not\in K$ with $a' \in K$, the resulting $H$-field $K(a)$ embeds over $K$ into any other $H$-field extension of $K$ containing an element $b$ with $b' = a'$ \emph{via an embedding sending $a$ to $b$}.
The italicized part is only true if we are in the setting of Fact~\ref{fact:largeint}; many thanks are due to Lou for noticing this gap.

Fixing the gap that this created in the proof of Proposition~\ref{prop:dtranselem} turned out to be quite involved.
Shortly after discovering the gap, I managed to prove this proposition in case~\ref{caseB} of Section~\ref{sec:vgext}, where $f>K$, modulo Proposition~\ref{prop:hensel}.
Near the end of 2024, I saw that this approach would also work for ``wide gaps'' in the sense of~\cite{ADH24C} and conjectured that $v(\log\fm_n-M_n) = (\lambda_n+\Z\lambda_n^\dagger)^\downarrow$ in the general case (see Section~\ref{sec:vgext} for notation).
I had a number of very helpful conversations with Vincent Bagayoko about this, and I thank Vincent for these conversations and for affording me the opportunity to visit Paris.

Despite this progress, things went no further until September 2026, when I gave the original preprint and all of the work I had done on closing this gap to Claude Opus 5.5 (provided through the Stockholm University mathematics department).
Claude managed to find a fix, and the main ideas in the proof of Proposition~\ref{prop:approxlog} come from this.
It is worth highlighting the most important insights to come out of this fix:
\begin{itemize}
\item The choice of compositional conjugation in Convention~\ref{conv:normalize}, which depends on a given $n$.
\item The derivations $\theta_m$, which provide useful decompositions of derivatives of elements in $K\langle f\rangle$.
\item The method of reducing approximations by elements in $K^{\h}$ to approximations by elements in $K$ via specialization in Proposition~\ref{prop:hensel}.
\item The key insight in Proposition~\ref{prop:twoints} that the rational Kaplansky Lemma of~\cite{Geh17C} can be played against logarithmic derivatives to show that adjoining one integral does not allow us to better approximate another.
\end{itemize}
Both Propositions~\ref{prop:hensel} and~\ref{prop:twoints} were originally proven in the specific setting of Section~\ref{sec:vgext}, though they may be of broader use and are included here as general results about $H$-fields.

The gap also rendered the proof of Proposition~\ref{prop:zgap} incomplete.
Claude noticed an error in an earlier attempted fix and suggested a shorter workaround, which constitutes the proof of the claim in this proposition.
The text of the article itself was written by me, and I take responsibility for its contents.
In addition to finding the proofs mentioned above, Claude Opus 5.5 was used for some initial restructuring and for proofreading.
I would like to thank Lou and Allen for feedback on the current draft.

\section{Preliminaries}\label{sec:prelim}
We draw heavily from~\cite{ADH17}.
The book~\cite{ADH17} is almost entirely self-contained, so some of the results we cite (especially the results on valued fields from~\cite[Chapter~3]{ADH17}) were not originally proven there.
We use the same notational conventions as~\cite{ADH17}, but we will repeat what is needed in this paper.

We let $k$, $m$, $n$ and $r$ range over $\N=\{0,1,2,\ldots\}$.
By ``ordered set'' we mean ``totally ordered set''.
Let $S$ be an ordered set.
For $a \in S$, we set $S^{>a}\coloneqq \{s\in S:s>a\}$; similarly for $S^{\geq a}$, $S^{<a}$, $S^{\leq a}$, and $S^{\neq a}$.
For $A \subseteq S$, we set
\[
A^{\downarrow} \coloneqq \{s\in S:s\leq a \text{ for some }a \in A\}
\]
and we say that \textbf{$A$ is downward closed in $S$} if $A = A^{\downarrow}$.

Now let $\Gamma$ be an ordered abelian group.
Then we let $\Gamma^{>}\coloneqq \Gamma^{>0}$ and we define $\Gamma^{\geq}$, $\Gamma^<$, $\Gamma^{\leq}$, and $\Gamma^{\neq}$ analogously.
For $\gamma \in \Gamma$, we let $[\gamma]$ denote the archimedean class of $\gamma$, and for $S\subseteq \Gamma$, we set $[S] \coloneqq\{[\gamma]:\gamma \in S\}$.
Given also $\alpha,\beta \in \Gamma$, we write $\alpha = \beta+o(\gamma)$ to mean that $\alpha = \beta+\delta$ for some $\delta \in \Gamma$ with $[\delta]<[\gamma]$; similarly for $\alpha\leq \beta+o(\gamma)$.
We consider $[\Gamma]$ as a totally ordered set with minimal element $[0]$.

Given a ring $R$, we let $R^\times$ denote the group of invertible elements in $R$.
All fields are assumed to be of characteristic zero.

\subsection{Ordered fields and ordered exponential fields}
Let $K$ be an ordered field.
We let $K^{\rc}$ denote the real closure of $K$.
If $L$ is a real closed ordered field extension of $K$, then there is a unique ordered field embedding $K^{\rc}\to L$ over $K$.

An \textbf{exponential} on $K$ is an ordered group isomorphism $\exp\colon K\to K^>$.
We define an \textbf{ordered exponential field} to be an ordered field $K$ equipped with an exponential $\exp$.
Given such a field, we use $\log\colon K^>\to K$ to denote the inverse of $\exp$.
If $K$ and $L$ are ordered exponential fields, then an ordered field embedding $\iota\colon K\to L$ is said to be an \textbf{ordered exponential field embedding} if $\iota(\exp a) = \exp \iota(a)$ for all $a \in K$.
We use $\R_{\exp}$ to denote the real exponential field.
Then $\Th(\R_{\exp})$, the first-order theory of $\R_{\exp}$ in the natural language of ordered exponential fields, is model complete and o-minimal by Wilkie's theorem~\cite{Wil96}.

\subsection{Valued fields}
Let $K = (K,\cO)$ be a valued field, so $\cO$ is a valuation ring with unique maximal ideal $\smallo$.
We let $\res(K) \coloneqq \cO/\smallo$ denote the residue field of $K$ and we let $\Gamma \coloneqq K^\times/\cO^\times$ denote the value group of $K$ with valuation $v \colon K^\times\to \Gamma$.
Our standing assumption that all fields are of characteristic zero applies to $\res(K)$ as well, so all valued fields are of equicharacteristic zero.
We construe $\Gamma$ as an additive group, totally ordered by $va> 0 \iff a\in \smallo$.
We extend $v$ to all of $K$ by setting $v(0) \coloneqq \infty>\Gamma$, and for $a,b \in K$ we set
\begin{align*}
a\asymp b &\logeq va = vb,& a\preceq b &\logeq va \geq vb,\\
a\prec b &\logeq va > vb,& a\sim b &\logeq a-b\prec a.
\end{align*}
Note that $a\sim b$ if and only if $a,b \neq 0$ and $a/b \in 1+\smallo$.

If $L$ is also a valued field, then we denote its value group, valuation ring, and maximal ideal by $\Gamma_L$, $\cO_L$, and $\smallo_L$ respectively.
Let $L$ be a valued field extension of $K$.
We identify $\Gamma$ with a subgroup of $\Gamma_L$ and we identify $\res(K)$ with a subfield of $\res(L)= \cO_L/\smallo_L$ in the natural way.
For $a \in L$, we let $v(a-K)$ denote the set $\{v(a-y):y \in K\} \subseteq \Gamma_L\cup\{\infty\}$.
We say that $L$ is an \textbf{immediate extension of $K$} if $\Gamma_L = \Gamma$ and $\res(L) = \res(K)$.

Let $(a_\rho)$ be a well-indexed sequence in $K$ (that is, a sequence indexed by ordinals below some infinite limit ordinal).
Then $(a_\rho)$ is a \textbf{pseudocauchy sequence} (or \textbf{pc-sequence}) if
\[
\tau>\sigma>\rho \implies v(a_\tau-a_\sigma)>v(a_\sigma-a_\rho)\quad\text{eventually},
\]
where \emph{eventually} means whenever the indices are greater than some $\rho_0$.
An element $a$ in a valued field extension of $K$ is a \textbf{pseudolimit} of $(a_\rho)$ if $v(a-a_\rho)$ is eventually strictly increasing; in this case, we write $a_\rho \leadsto a$.
If $P(a_\rho) \not\leadsto 0$ for all $P \in K[X]\setminus K$, then $(a_\rho)$ is said to be of \textbf{transcendental type} over $K$.
A \textbf{divergent} pc-sequence in $K$ is one with no pseudolimit in $K$; so any pc-sequence of transcendental type is necessarily divergent.
The following is the ``rational Kaplansky Lemma'' of~\cite{Geh17C}.

\begin{fact}[{\cite[Lemma~4.4]{Geh17C}}]\label{fact:kaplansky}
Let $L$ be a valued field extension of $K$.
Let $a \in L\setminus K$ and let $(a_\rho)$ be a divergent pc-sequence in $K$ of transcendental type with $a_\rho\leadsto a$.
Let $F \in K(X)\setminus K$, and put $\gamma_\rho\coloneqq v(a_\rho-a)$.
Then the following hold eventually (i.e.\ for all $\rho$ greater than some $\rho_0$):
\begin{enumerate}[label=\normalfont(\roman*), align = right]
\item\label{kaplanskyi} $F(a_\rho)$ is a well-defined element of $K$ (i.e.\ the denominator is nonzero).
\item\label{kaplanskyii} $(F(a_\rho))$ is a divergent pc-sequence in $K$, and $F(a_\rho) \leadsto F(a)$.
Thus $v F(a_\rho) = vF(a)$.
\item\label{kaplanskyiii} $v(F(a_\rho)- F(a)) = \alpha + k\gamma_\rho$ for some $\alpha \in \Gamma_L$ and some $k \geq 1$.
\end{enumerate}
\end{fact}
In~\cite{Geh17C}, $L$ is assumed to be an immediate extension of $K$, but the proofs of the parts listed here do not use this.
The following improvement on Fact~\ref{fact:kaplansky}\ref{kaplanskyiii} is an unpublished lemma of Gehret; we give a short proof:

\begin{lemma}\label{lem:mvt}
In the last part of Fact~\ref{fact:kaplansky}, we can always take $\alpha = vF'(a)$ and $k = 1$.
\end{lemma}
\begin{proof}
Write $F = P/Q$ with $P,Q\in K[X]$ and set $R\coloneqq P-F(a)Q\in L[X]$, so
\[
R(a) = 0,\qquad R' = QF'+(F-F(a))Q',\qquad F-F(a) = R/Q.
\]
Taking Taylor expansions of $R$ and $R'$ around $a$, we have
\[
R(a_\rho) = \sum_{i\geq1}\frac{1}{i!}R^{(i)}(a)(a_\rho-a)^i,\qquad R'(a_\rho) = \sum_{i\geq1}\frac{1}{(i-1)!}R^{(i)}(a)(a_\rho-a)^{i-1}.
\]
By~\cite[Lemma~3.2.3]{ADH17} (and the fact that $v(1/i!) = 0$ in equicharacteristic zero), we find $k\geq 1$ such that eventually,
\[
v(R(a_\rho)) = v(R^{(k)}(a))+k\gamma_\rho<v(R^{(i)}(a))+i\gamma_\rho\ \text{ for $i\neq k$}.
\]
For brevity, we say that $v(R(a_\rho))$ is eventually \emph{linear with slope $k$}.
Note that by the Taylor expansion of $R'$ and the inequality above, $v(R'(a_\rho))$ is eventually linear with slope $k- 1$.
On the other hand
\[
R'(a_\rho) = Q(a_\rho)F'(a_\rho)+(F(a_\rho)-F(a))Q'(a_\rho),
\]
where the valuation of the first term is eventually constant by Fact~\ref{fact:kaplansky}\ref{kaplanskyii} applied to $Q$ and $F'$ and the valuation of the second is either infinite or eventually linear with slope $k$.
As $(\gamma_\rho)$ is eventually strictly increasing, we see that $v(R'(a_\rho))$ must either be eventually constant or linear with slope $k$.
This is only possible if $k-1 = 0$, so $k = 1$.
Thus eventually,
\[
v(F(a_\rho)-F(a)) = v(R(a_\rho))- v(Q(a_\rho)) = v R'(a) + \gamma_\rho- vQ(a) = vF'(a)+ \gamma_\rho.\qedhere
\]
\end{proof}

We let $K^{\h}$ denote the henselization of $K$.
Then $K^{\h}$ is an immediate extension of $K$, and there is a unique valued field embedding of $K^{\h}$ over $K$ into any henselian valued field extending $K$~\cite[Proposition~3.3.25 and Corollary~3.3.26]{ADH17}.
Accordingly, if $L$ is a henselian valued field extension of $K$, we identify $K^{\h}$ with a valued subfield of $L$.
As $K$ has equicharacteristic zero, any divergent pc-sequence in $K^{\h}$ is necessarily of transcendental type~\cite[Corollary~3.3.21]{ADH17}.
A key fact for us is Kuhlmann's approximation theorem for elements in $K^{\h}$:

\begin{fact}[{\cite[Theorem~1]{Kuh11}}]\label{fact:kuhlmann}
Let $x \in K^{\h}\setminus K$.
Then there is a nontrivial convex subgroup $\Delta \subseteq \Gamma$ and an element $\alpha \in \Gamma$ with
\[
v(x-K) = (\alpha+\Delta)^{\downarrow}.
\]
\end{fact}

Finally, we need the following fact about ordered valued fields:

\begin{fact}[{\cite[Corollary~3.5.19]{ADH17}}]\label{fact:realclosed}
Suppose that $K$ is an ordered field and that $\cO$ is convex with respect to the ordering on $K$.
Then $K$ is real closed if and only if $K$ is henselian, $\res(K)$ is real closed, and $\Gamma$ is divisible.
\end{fact}

\subsection{Differential fields}
Let $K$ be a differential field with derivation $\der\colon K \to K$.
We let $K\{Y\}$ denote the ring of differential polynomials over $K$.
Then $K\{Y\}$ is a differential ring extension of $K$.
For $a \in K$ we often write $a'$ instead of $\der(a)$ and $a''$ instead of $\der^2(a)$.
More generally, we write $a^{(n)}$ for $\der^n(a)$.
We set $a^\dagger \coloneqq a'/a$ if $a \neq 0$ and we call $a^\dagger$ the \textbf{logarithmic derivative} of $a$.

We let $C\coloneqq \ker(\der)$ denote the constant field of $K$ and if $L$ is also a differential field, then we denote its constant field by $C_L$.
Given $\phi \in K^\times$, we let $K^\phi$ be the differential field with underlying field $K$ and derivation $\phi\inv\der$.
We call $K^{\phi}$ the \textbf{compositional conjugate of $K$ by $\phi$}.
Note that $C_{K^\phi} = C$.

We set $K' \coloneqq \{a':a \in K\}$ and we set $(K^\times)^\dagger\coloneqq \{a^\dagger:a \in K^\times\}$.
Note that $(K^\times)^\dagger$ is an additive subgroup of $K$ since $a^\dagger + b^\dagger = (ab)^\dagger$ for $a,b \in K^\times$.
Note also that $K'$ is a $C$-vector subspace of $K$.
We say that $a\in K$ can be \textbf{integrated} if $a \in K'$ and we call an element $b \in K$ with $b' = a$ an \textbf{integral} of $a$.
If $b_0$ and $b_1$ are integrals of $a$, then $b_0-b_1 \in C$.
We say that $a\in K$ can be \textbf{exponentially integrated} if $a \in (K^\times)^\dagger$ and we call an element $b \in K^\times$ with $b^\dagger = a$ an \textbf{exponential integral} of $a$.
If $b_0$ and $b_1$ are exponential integrals of $a$, then $b_0/b_1 \in C^\times$.

For $y$ in a differential field extension of $K$, we let $K\langle y\rangle\coloneqq K(y,y',y'',\ldots)$ denote the differential field extension of $K$ generated by $y$.
We say that $y$ is \textbf{$\d$-transcendental over $K$} if the sequence $y,y',y'',\ldots$ is algebraically independent over $K$ (in the field-theoretic sense) and we say that $y$ is \textbf{$\d$-algebraic over $K$} otherwise.
Note that $y$ is $\d$-algebraic over $K$ if and only if $K\langle y\rangle$ has finite transcendence degree over $K$.

For $P \in K[X]$, we let $P^\der$ denote the polynomial obtained by taking the derivatives of the coefficients of $P$.
Then $P(y)' = P^\der(y)+P'(y)y'$ for any $y$ in a differential field extension of $K$.
We extend this to rational functions $F \in K(X)$ using the quotient rule; then $F^\der$ still has the property $F(y)' = F^\der(y)+F'(y)y'$.
The following fact is proven just before~\cite[Corollary~5.1.28]{ADH17}.

\begin{fact}\label{fact:trace}
Let $A(Y) = g$ be a linear differential equation over $K$.
If this equation has a solution in a differential field extension of $K$ that is algebraic over $K$, then it has a solution in $K$.
\end{fact}

\section{\texorpdfstring{$H$}{H}-fields}\label{sec:asymp}
Let $K$ be an ordered differential field with constant field $C$.
Following~\cite{AD02}, we say that $K$ is an \textbf{$H$-field} if
\begin{enumerate}[label=\normalfont(H\arabic*), align = right]
\item\label{H1} $f'>0$ for all $f \in K$ with $f>C$;
\item\label{H2} $\cO = C+\smallo$, where $\cO$ is the convex hull of $C$ in $K$ and $\smallo$ is the unique maximal ideal of $\cO$.
\end{enumerate}
We always view an $H$-field as an ordered \emph{valued} differential field with the convex valuation ring $\cO$ as in~\ref{H2}.
An ordered field with trivial derivation is an $H$-field with trivial valuation, and if the derivation of an $H$-field is nontrivial, then so is its valuation.

In this section, we recall some needed facts on $H$-fields, asymptotic couples, and newtonizations from~\cite{ADH17}.
We also prove in Propositions~\ref{prop:hensel} and~\ref{prop:twoints} below that neither passing to the henselization nor adjoining a suitably different integral improves approximations to an integral in an $H$-field.

\begin{convention*}
For the remainder of this paper, $K$ is an $H$-field with valuation ring $\cO$, maximal ideal $\smallo$, derivation $\der$, constant field $C$, and value group $\Gamma$.
\end{convention*}

By~\ref{H2}, the residue map sends $C$ isomorphically onto $\res(K)$.
Thus an $H$-field extension $L$ of $K$ is immediate if and only if $\Gamma_L = \Gamma$ and $C_L = C$.
By~\cite[Proposition~9.5.3 and Lemma~10.5.8]{ADH17}, we always view the henselization of an $H$-field $K$ as an immediate $H$-field extension of $K$.
Also by~\ref{H2}, each $u \in \cO$ differs from an element in $\smallo$ by a constant, so $\der \cO = \der \smallo$.
We say that $K$ has \textbf{small derivation} if $\der\smallo\subseteq\smallo$, and \textbf{large derivation} otherwise.
Note that if $\phi\in K^>$, then the compositional conjugate $K^\phi$ is again an $H$-field.

\subsection{Asymptotic couples}\label{subsec:Hfields}
By~\cite[Lemma~10.5.2]{ADH17}, the $H$-field $K$ is \textbf{$H$-asymptotic}:\ for all nonzero $f,g \in \smallo$, we have
\[
f\prec g \Longleftrightarrow f' \prec g',\qquad f\preceq g\Longrightarrow f^\dagger \succeq g^\dagger.
\]
Thus, for $f \in K^\times \setminus \cO^\times$, the values $v(f')$ and $v(f^\dagger)$ depend only on $vf$.
For $\gamma = vf$ with $\gamma\neq 0$, we write
\[
\gamma^\dagger\coloneqq v(f^\dagger),\qquad \gamma'\coloneqq v(f') = \gamma+\gamma^\dagger.
\]
Following Rosenlicht~\cite{Ros81}, we denote the map $\gamma\mapsto\gamma^\dagger\colon \Gamma^{\neq}\to\Gamma$ by $\psi$, and we call the pair $(\Gamma,\psi)$ the \textbf{asymptotic couple of $K$}.

\begin{fact}[{\cite[Lemma~6.5.4(ii),(iii)]{ADH17}}]\label{fact:6.5.4}
$(\Gamma,\psi)$ has the following properties:
\begin{enumerate}[label=\normalfont(\roman*), align = right]
\item\label{6.5.4i} If $\alpha,\beta \in \Gamma^{\neq}$ and $\alpha \neq \beta$, then $\psi(\alpha)-\psi(\beta) = o(\alpha-\beta)$.
\item\label{6.5.4ii} The map $\gamma\mapsto \gamma'\colon \Gamma^{\neq}\to \Gamma$ is strictly increasing.
\end{enumerate}
\end{fact}

We set
\[
(\Gamma^<)' \coloneqq \{\gamma':\gamma \in \Gamma^<\},\qquad (\Gamma^>)' \coloneqq \{\gamma':\gamma \in \Gamma^>\},\qquad \Psi \coloneqq \psi(\Gamma^{\neq}).
\]
Then $(\Gamma^<)'$ is downward closed, $(\Gamma^>)'$ is upward closed, and $(\Gamma^<)'<(\Gamma^>)'$.
Note that $K$ has small derivation if and only if $(\Gamma^>)' \subseteq \Gamma^>$, and large derivation if and only if $0 \in (\Gamma^>)'$.
We say that $K$ has \textbf{asymptotic integration} if $\Gamma = (\Gamma^<)' \cup (\Gamma^>)'$.
If $K$ does not have asymptotic integration, then there is only one $\beta \in \Gamma$ with $(\Gamma^<)'<\beta <(\Gamma^>)'$.
In this case, either $\beta = \max \Psi$ and $K$ is said to be \textbf{grounded}, or $\beta \not\in \Psi$ but $\beta = \sup \Psi$ and we call $\beta$ a \textbf{gap in $K$}.
If $K$ has trivial valuation, then $\Psi = (\Gamma^<)' = (\Gamma^>)' = \emptyset$, and $0$ is a gap in $K$.

\begin{fact}[{\cite[Corollary~9.2.16]{ADH17}}]\label{fact:trich}
Exactly one of the following holds:
\begin{enumerate}[label=\normalfont(\roman*), align = right]
\item\label{trichi} $K$ has asymptotic integration;
\item\label{trichii} $K$ has a gap;
\item\label{trichiii} $K$ is grounded.
\end{enumerate}
\end{fact}

Suppose $K$ is ungrounded.
Then $\Psi \subseteq (\Gamma^{<})'$, so for $\gamma \in \Gamma^<$, there is a unique element $\chi(\gamma)\in \Gamma^<$ with $\chi(\gamma)'= \gamma^\dagger$.
This map $\chi\colon \Gamma^<\to \Gamma^<$ is called the \textbf{contraction map}.

\subsection{Integrals in \texorpdfstring{$H$}{H}-fields}\label{subsec:inthensel}
The following result plays a key role in the proof of Proposition~\ref{prop:approxlog}.
It tells us that integrals of elements in $K$ cannot be better approximated by passing to the henselization of $K$.

\begin{proposition}\label{prop:hensel}
Let $f$ be an element in an $H$-field extension of $K$ with $f' \in K$.
Then
\[
v(f-K^{\h}) = v(f-K).
\]
\end{proposition}
\begin{proof}
This is trivial if $f \in K$, so suppose $f \not\in K$.
We work in a henselian $H$-field extension $M$ of $K$ containing $f$, and we identify $K^{\h}$ with an $H$-subfield of $M$.
Set $S\coloneqq v(f-K)$, let $y\in K^{\h}$, and suppose towards a contradiction that $v(f-y)>S$.
Then, $v(y-K) = S$, so Fact~\ref{fact:kuhlmann} yields $\alpha\in\Gamma$ and a nontrivial convex subgroup $\Delta\subseteq\Gamma$ with $S = (\alpha+\Delta)^{\downarrow}$.
We may assume $\alpha \neq 0$, and we take $a,e \in K$ with $va = v(f-e) = \alpha$.
Set
\[
w \coloneqq a\inv(f-e),\qquad x \coloneqq a\inv(y-e) \in K^{\h},\qquad \phi\coloneqq a^\dagger\in K,\qquad h\coloneqq (f-e)'/a'\in K.
\]
Then $vw = vh = 0$ and $v(w-x) = v(f-y)-\alpha>\Delta$.
As $\alpha\neq0$, we have $v\phi\in\Psi$.
For the remainder of the proof, we work in $M^\phi$, so we have small derivation, $a^\dagger = 1$, and $w' + w = h$.
Note that $\phi$ may be negative, so $M^\phi$ need not be an $H$-field, but it is still \emph{$\d$-valued} with the same valuation ring and constant field as $M$ (see~\cite[Section~9.1]{ADH17}), and this is all that we use below.
We claim that
\begin{equation}\label{eq:invariant}
(S\setminus \{0\})' \subseteq S.
\end{equation}
Since the map $\gamma\mapsto\gamma'$ on $\Gamma_M^{\neq}$ is strictly increasing by Fact~\ref{fact:6.5.4}\ref{6.5.4ii}, it suffices to show that $(\alpha+\delta)' \in S$ for $\delta \in \Delta^{\neq}$ with $\alpha +\delta\neq 0$.
For this, we use that $\psi_M(\alpha) = 0$, together with Fact~\ref{fact:6.5.4}\ref{6.5.4i} to get
\[
(\alpha+ \delta)' = \alpha + \delta + \psi_M(\alpha+\delta) = \alpha + \delta + \psi_M(\alpha+\delta)- \psi_M(\alpha) = \alpha + \delta+ o(\delta) \in \alpha+\Delta.
\]

Now let $\dot\cO\coloneqq \{u\in M:vu\geq\delta\text{ for some }\delta\in\Delta\}$, so $\dot\cO$ is a convex valuation ring of $M$ with maximal ideal $\dot\smallo\coloneqq \{u \in M:vu>\Delta\}$.
Let $\dot M$ denote the residue field $\dot \cO/\dot\smallo$, so the derivation on $M$ induces a derivation on $\dot M$ by~\cite[Lemma~4.4.2 and Corollary~4.4.4]{ADH17}.
We let $\dot K\subseteq\dot M$ be the differential residue field of $K\cap\dot\cO$.
The elements $w$, $x$, and $h$ all lie in $\dot\cO$, and we let $\dot w$, $\dot x$, and $\dot h$ denote their images in $\dot M$, so $\dot w = \dot x$ and $\dot w' + \dot w = \dot h \in \dot K$.
As $x$ is algebraic over $K$, its residue $\dot w = \dot x$ is algebraic over $\dot K$, so Fact~\ref{fact:trace} provides $\dot z\in\dot K$ with $\dot z'+\dot z = \dot h$.
Lift $\dot z$ to an element $z\in K$, and put $g \coloneqq f-e-az \in f-K$, so $vg \in S$.
As we have arranged that $a' = a$, we see that
\begin{equation}\label{eq:gprime}
g' = f'-e'-a'z-az' = ah-az-az' \in a\dot\smallo.
\end{equation}
Using that $S$ has no greatest element, we arrange that $vg\in S\setminus \{0\}$ by subtracting a constant $c \in C$ with $g\sim c$ if need be; note that this doesn't change $g'$.
Then~\eqref{eq:invariant} gives $vg'\in S$, contradicting~\eqref{eq:gprime}.
\end{proof}

The next proposition shows that if two integrals of elements in $K$ differ (as measured by the eventual value of $\psi$ on $v(a-K)$ and $v(b-K)$), then adjoining one does not allow us to better approximate the other.

\begin{proposition}\label{prop:twoints}
Let $a,b$ be elements in an $H$-field extension of $K$ with $a,b\not\in K$ and $a',b' \in K$.
Put
\[
S_a\coloneqq v(a-K),\qquad S_b\coloneqq v(b-K).
\]
Suppose that $S_a$ has no largest element and that there are $\alpha\neq\beta$ such that $\gamma^\dagger = \alpha$ for all sufficiently large $\gamma\in S_a$ and $\gamma^\dagger = \beta$ for all sufficiently large $\gamma\in S_b$.
Then $v(b-K(a)) = S_b$.
\end{proposition}
\begin{proof}
Using that $S_a$ has no largest element, take a sequence $(a_\rho)$ in $K$ such that $\gamma_\rho\coloneqq v(a-a_\rho)$ is strictly increasing and cofinal in $S_a$.
Then $(a_\rho)$ is a divergent pc-sequence in $K$ and $a_\rho\leadsto a$.
By Proposition~\ref{prop:hensel}, $(a_\rho)$ remains divergent in $K^{\h}$, so $(a_\rho)$ is of transcendental type and $K(a)$ is an immediate extension of $K$.
Suppose towards a contradiction that there is $y\in K(a)$ with $v(b-y)> S_b$ (we allow for the possibility that $y = b$, so $v(b-y)$ may be infinite).
Write $y = F(a)$ for $F \in K(X)\setminus K$.
Passing to an end segment, we may assume that $F(a_\rho)$ is defined for each $\rho$, and we set
\[
y_\rho\coloneqq F(a_\rho) \in K,\qquad \eta_\rho\coloneqq v(y_\rho-y),\qquad \delta \coloneqq vF'(a).
\]
Then $(y_\rho)$ is a divergent pc-sequence in $K$ with $y_\rho \leadsto y$ by Fact~\ref{fact:kaplansky}\ref{kaplanskyii}.
After possibly passing to an end segment again, Lemma~\ref{lem:mvt} gives
\begin{equation}\label{eq:mvt}
\eta_\rho = \delta + \gamma_\rho.
\end{equation}
As $v(b-y_\rho)\in S_b$ and $v(b-y)>S_b$, we see that $v(b-y_\rho) = v(y-y_\rho) = \eta_\rho$,
so $y_\rho \leadsto b$.
As $(y_\rho)$ is divergent, $(\eta_\rho)$ is cofinal in $S_b$.
Passing to an end segment, we may arrange that
\[
\gamma_\rho\neq 0,\qquad \eta_\rho\neq 0,\qquad \gamma_\rho^\dagger = \alpha,\qquad \eta_\rho^\dagger = \beta
\]
for all $\rho$.

Set $G\coloneqq F^\der+ a' F' \in K(X)$, so $y' = F(a)' = G(a)$.
Then
\begin{equation}\label{eq:compare}
y_\rho' - G(a_\rho) = F'(a_\rho)(a_\rho'-a').
\end{equation}
We now calculate the valuations on both sides of~\eqref{eq:compare}.
Applying Fact~\ref{fact:kaplansky}\ref{kaplanskyii} with $F'$ in place of $F$ and using $\gamma_\rho\neq 0$, we have
\[
v\big(F'(a_\rho)(a_\rho'-a')\big) = \delta+ \gamma_\rho' = \delta+\gamma_\rho+\alpha.
\]
Note that Fact~\ref{fact:kaplansky} only applies to nonconstant rational functions, but the equality above is trivial if $F' \in K$.
For the left side, we have $y_\rho' - G(a_\rho) = (y_\rho' -b') - (G(a_\rho)- b')$, and we note that
\[
v(y_\rho' -b') = \eta_\rho',\qquad v(G(a_\rho)-b') = v(G(a)- b') = v(y'- b') > \eta_\rho',
\]
where the equality on the right uses Fact~\ref{fact:kaplansky}\ref{kaplanskyii} with $G- b'$ in place of $F$ (again, this is trivial if $G-b'\in K$) and the inequality uses that $v(y - b)>\eta_\rho$, as $\eta_\rho \in S_b$.
Thus, $v(y_\rho' - G(a_\rho))= \eta_\rho' = \eta_\rho+\beta$, and~\eqref{eq:compare} yields
\[
\eta_\rho+\beta = \delta+ \gamma_\rho+\alpha.
\]
As $\alpha\neq \beta$, this contradicts~\eqref{eq:mvt}.
\end{proof}

\begin{remark}\label{rem:dvalued}
Propositions~\ref{prop:hensel} and~\ref{prop:twoints} hold in greater generality:\ we only need that $K$ is a \emph{$\d$-valued field} in the sense of~\cite[Section~9.1]{ADH17} and that $f,a,b$ lie in a $\d$-valued field extension of $K$.
The same proofs work verbatim.
\end{remark}

We will need the following two facts about adjoining integrals to henselian $H$-fields.

\begin{fact}[{\cite[Lemmas~10.2.4,~10.2.5(iii), and~10.5.8]{ADH17}}]\label{fact:smallint}
Let $K$ be henselian and let $s \in K$ with $vs \in (\Gamma^>)'$ and $s\notin \der\smallo$.
Then there is an immediate $H$-field extension $K(a)$ of $K$ with $a\prec 1$ and $a' = s$.
If $L$ is also an $H$-field extension of $K$ with $s \in \der\smallo_L$, then there is a unique $H$-field embedding $K(a)\to L$ over $K$.
\end{fact}

\begin{fact}[{\cite[Lemmas~10.2.5(ii),~10.2.6, and~10.5.8]{ADH17}}]\label{fact:largeint}
Let $K$ be henselian with asymptotic integration and let $s \in K\setminus K'$ with $v(s-y')<(\Gamma^>)'$ for all $y \in K$.
Then there is an immediate $H$-field extension $K(a)$ of $K$ with $a\succ 1$ and $a' = s$.
If $L$ is also an $H$-field extension of $K$ and $a^* \in L$ satisfies $(a^*)' =s$, then there is a unique $H$-field embedding $K(a)\to L$ over $K$ that sends $a$ to $a^*$.
\end{fact}

Note that in Facts~\ref{fact:smallint} and~\ref{fact:largeint}, the element $a$ is necessarily transcendental over $K$ by Fact~\ref{fact:trace}.
Iterating these facts along with taking henselizations, we have the following:

\begin{fact}[{\cite[Proposition~10.2.7 and Lemma~10.5.8]{ADH17}}]\label{fact:intclosure}
Suppose $K$ has asymptotic integration.
Then $K$ has
an immediate $H$-field extension $K(\int)$ that is henselian, closed under integration, and embeds over $K$ into any henselian $H$-field extension of $K$ that is closed under integration.
\end{fact}

Finally, we show that any $H$-field has a canonical immediate extension that admits a definable logarithm $1+\smallo\to \smallo$.

\begin{lemma}\label{lem:ellext}
$K$ has an immediate $H$-field extension $K^\ell$ such that
\begin{enumerate}[label=\normalfont(\roman*), align = right]
\item\label{ellexti} $K^\ell$ is henselian;
\item\label{ellextii} $(1+\smallo_{K^\ell})^\dagger\subseteq \der\smallo_{K^\ell}$;
\item\label{ellextiii} $K^\ell$ embeds uniquely over $K$ into any henselian $H$-field extension $L$ of $K$ with $(1 +\smallo_L)^\dagger\subseteq\der\smallo_L$.
\end{enumerate}
\end{lemma}
\begin{proof}
We define an \emph{$\ell$-tower} on $K$ to be an increasing chain of henselian $H$-fields $(K_\mu)_{\mu\leq\nu}$ such that $K_0 = K^{\h}$, $K_\mu = \bigcup_{\lambda<\mu}K_\lambda$ for each limit ordinal $\mu$, and $K_{\mu+1} = K_\mu(y_\mu)^{\h}$ where $y_\mu\prec 1$ and $y_\mu' = a_\mu^\dagger$ for some $a_\mu\in 1+\smallo_{K_\mu}$ with $a_\mu^\dagger\notin \der\smallo_{K_\mu}$.
Note that if $(K_\mu)_{\mu\leq\nu}$ is an $\ell$-tower on $K$, then each $K_\mu$ is an immediate extension of $K$ by Fact~\ref{fact:smallint}.
Thus, maximal $\ell$-towers on $K$ exist, and we let $(K_\mu)_{\mu\leq\nu}$ be such a tower.
Then $K_\nu$ is a henselian immediate extension of $K$ and $(1+\smallo_{K_\nu})^\dagger\subseteq \der\smallo_{K_\nu}$.
If $L$ is a henselian $H$-field extension of $K$ with $(1 +\smallo_L)^\dagger\subseteq\der\smallo_L$, then the universal property of the henselization, along with the universal property in Fact~\ref{fact:smallint}, gives a unique $H$-field embedding $K_\mu\to L$ over $K$ for each $\mu\leq \nu$.
Thus, we may take $K^\ell\coloneqq K_\nu$.
\end{proof}

Suppose $K = K^\ell$.
If $E$ is an $H$-subfield of $K$, then we identify $E^\ell$ with its unique image in $K$.
We define a function $\ln\colon 1+\smallo\to \smallo$ by letting $\ln a$ be the unique element of $\smallo$ satisfying
\[
(\ln a)' = a^\dagger
\]
for each $a \in 1+\smallo$.
Given $\epsilon \in \smallo^{\neq}$, we have $\ln(1+\epsilon)' \sim \epsilon'$, so $\epsilon >0$ if and only if $\ln(1+\epsilon) > 0$.
From this it is straightforward to check that $\ln$ is an ordered group embedding $1+\smallo\to \smallo$.

\subsection{\texorpdfstring{$\upo$}{ω}-free \texorpdfstring{$H$}{H}-fields and newtonianity}\label{subsec:upofree}
As shown in~\cite[Corollary~16.2.5]{ADH17}, the theory of $H$-fields has a model companion, namely the theory $T^{\nl}$ of $\upo$-free newtonian Liouville closed $H$-fields.
We discuss these axioms below, but first we note that they are invariant under compositional conjugation:\ if $K \models T^{\nl}$, then $K^\phi\models T^{\nl}$ for any $\phi>0$.

We say that $K$ is \textbf{Liouville closed} if $K$ is real closed and every $a \in K$ can be integrated and exponentially integrated in $K$.
The axiom of \textbf{$\upo$-freeness} is more technical; fortunately, this paper only uses its consequences, which we describe below.
If $K$ is $\upo$-free, then it is ungrounded by definition and it does not have a gap~\cite[Corollaries~11.6.1 and~11.7.3]{ADH17}, so $\upo$-free $H$-fields have asymptotic integration.
The utility of $\upo$-freeness is that it is quite robust:\ it transfers to $\d$-algebraic extensions and to certain subfields:

\begin{fact}[{\cite[Section~11.7 and Theorem~13.6.1]{ADH17}}]\label{fact:dalgupofree}
Let $K$ be $\upo$-free, let $L$ be a $\d$-algebraic $H$-field extension of $K$, and let $E$ be an $H$-subfield of $K$ with $\Gamma_E^<$ cofinal in $\Gamma^<$.
Then $L$ and $E$ are $\upo$-free as well.
\end{fact}

The axiom of \textbf{newtonianity} is also quite technical, and as with $\upo$-freeness, we will only use its consequences.
One may think of it as a differential analog of henselianity (more precisely, an approximate version of differential henselianity).
When coupled with $\upo$-freeness, it is quite a strong axiom:

\begin{fact}\label{fact:asdmax}
If $K$ is $\upo$-free and newtonian, then $K$ is \textbf{asymptotically $\d$-algebraically maximal}, that is, $K$ has no proper immediate $\d$-algebraic $H$-field extension.
\end{fact}

Fact~\ref{fact:asdmax} was shown under the assumption that $K$ also has divisible value group~\cite[Theorem~14.0.2]{ADH17}, but this divisibility assumption can be removed; see~\cite{PC19}.
Here is a consequence of Facts~\ref{fact:asdmax} and~\ref{fact:intclosure} that we use repeatedly:

\begin{fact}[{\cite[Lemma~1.6]{Kap26}}]\label{fact:newtonint}
If $K$ is $\upo$-free and newtonian, then $K$ is henselian and every element of $K$ has an integral in $K$.
In particular, $K = K^\ell$.
\end{fact}

Combining Fact~\ref{fact:asdmax} with~\cite[Lemmas~2.4.2, 10.5.8, 11.4.7, and~11.4.8]{ADH17} gives us the following:

\begin{fact}\label{fact:immextmap}
Let $K$ be $\upo$-free and newtonian and let $y$ be an element in an $H$-field extension of $K$ such that $v(y-K)$ has no largest element.
Then $K\langle y \rangle$ is an immediate extension of $K$.
If $L$ is also an $H$-field extension of $K$ and $y^* \in L$ realizes the same cut as $y$ over $K$, then there is an $H$-field embedding $K\langle y \rangle \to L$ over $K$ that sends $y$ to $y^*$.
\end{fact}

If $K$ is $\upo$-free, then a \textbf{newtonization of $K$} is by definition an immediate newtonian $H$-field extension of $K$ that embeds over $K$ into any newtonian $H$-field extension of $K$.

\begin{fact}\label{fact:ntexist}
Suppose $K$ is $\upo$-free.
Then $K$ has a newtonization $K^{\nt}$ that is $\d$-algebraic over $K$, hence also $\upo$-free by Fact~\ref{fact:dalgupofree}.
Any two newtonizations of $K$ are isomorphic over $K$.
\end{fact}

Fact~\ref{fact:ntexist} was shown under the assumption that $K$ also has divisible value group~\cite[Corollaries~14.3.12 and~14.5.4]{ADH17}, but again, this divisibility assumption can be removed by~\cite{PC19}.

If $K$ is $\upo$-free, then a \textbf{Newton-Liouville closure of $K$} is by definition a newtonian Liouville closed $H$-field extension of $K$ that embeds over $K$ into any newtonian Liouville closed $H$-field extension of $K$.

\begin{fact}[{\cite[Corollaries~14.5.10 and~16.2.2]{ADH17}}]\label{fact:nlexist}
Suppose that $K$ is $\upo$-free.
Then $K$ has a Newton-Liouville closure $K^{\nl}$ that is $\d$-algebraic over $K$.
The constant field of $K^{\nl}$ is a real closure of $C$, and any other Newton-Liouville closure of $K$ is isomorphic to $K^{\nl}$ over $K$.
\end{fact}

As shown in~\cite[Theorem~16.0.3]{ADH17}, $\upo$-free newtonian Liouville closed $H$-fields have no proper $\d$-algebraic $H$-field extensions with the same constant field.
We can use this to give another characterization of the Newton-Liouville closure of $K$.

\begin{fact}[{\cite[Lemma~1.11]{Kap26}}]\label{fact:altnl}
Let $K$ be $\upo$-free and let $L$ be a newtonian Liouville closed $H$-field extension of $K$.
Then $L$ is a Newton-Liouville closure of $K$ if and only if $L$ is $\d$-algebraic over $K$ and $C_L$ is algebraic over $C$.
\end{fact}

\section{Logarithmic and exponential \texorpdfstring{$H$}{H}-fields}\label{sec:logexpfields}
In this section, we investigate $H$-fields that are equipped with a logarithm.
Logarithms on $H$-fields were previously studied in~\cite{AD05}.

\begin{definition}\label{def:logarithm}
A \textbf{logarithm on $K$} is an ordered group embedding $\log\colon K^{>}\to K$ such that
\begin{enumerate}[label=\normalfont(L\arabic*), align = right]
\item\label{L1} $\log(1+\smallo)\subseteq \smallo$;
\item\label{L2} $(\log f)'=f^\dagger$ for all $f\in K^>$.
\end{enumerate}
\end{definition}

We use $\exp$ to denote the inverse of $\log$ where it is defined.
Let $\log\colon K^{>}\to K$ be any map satisfying~\ref{L2}, and let $f \in K^>$.
If $f \in C^>$, then $(\log f)' = f^\dagger = 0$, so $\log$ restricts to a map $C^>\to C$.
We call this map the \textbf{induced logarithm on $C$}.

\begin{lemma}\label{lem:logasymp}
Let $\log\colon K^>\to K$ be a group morphism satisfying~\ref{L1} and~\ref{L2}, and let $f\in K^>$.
Then $f\asymp 1$ if and only if $\log f\preceq 1$.
\end{lemma}
\begin{proof}
If $f\asymp 1$, then $f = c(1+\epsilon)$ with $c\in C^>$ and $\epsilon\in\smallo$, so $\log f = \log c+\log(1+\epsilon)\preceq 1$ by~\ref{L1}.
Conversely, if $\log f \in \cO$, then $v(f^\dagger) = v(\log(f)') \in (\Gamma^>)'\cup\{\infty\}>\Psi$, so $vf = 0$.
\end{proof}

Let $\log$ be a logarithm on $K$.
If $f\in K^>$ and $f \succ 1$, then $\log f\succ 1$ by Lemma~\ref{lem:logasymp}, so $\psi(vf) =v(\log f)' \in (\Gamma^<)'$; in particular, $K$ is ungrounded.
The trace of the logarithm on $\Gamma^<$ is given by the contraction map:\ $v(\log y) = \chi(vy)$ for $y \in K^>$ with $vy<0$.
By~\cite[Lemma~9.2.18(iii)]{ADH17} we have $\gamma<n\chi(\gamma)<0$ for each $\gamma \in \Gamma^<$ and each $n>0$, so $\log$ is sufficiently ``slow-growing'':
\[
y\succ (\log y)^n\succ 1 \text{ for all }y \in K^>\text{ with }y \succ 1.
\]

\begin{lemma}\label{lem:autoorder}
Let $\log\colon K^>\to K$ be a group morphism that satisfies~\ref{L1} and~\ref{L2}.
Suppose that $\log|_{C^>}\colon C^>\to C$ is an ordered group embedding.
Then $\log$ is a logarithm on $K$.
\end{lemma}
\begin{proof}
Let $f\in K$ with $f>1$; we need to show that $\log f>0$.
First, consider the case that $f \succ 1$.
Then $\log f\succ 1$ by Lemma~\ref{lem:logasymp} and $(\log f)' = f^\dagger > 0$ by~\ref{H1} and~\ref{L2}, so $\log f$ is positive, again by~\ref{H1}.
Now suppose $f \asymp 1$, so $f = c(1+\epsilon)$ for some $c \in C^{\geq 1}$ and some $\epsilon \in \smallo$ and
\[
\log f = \log c + \log(1+\epsilon).
\]
If $c>1$, then $\log c > 0$ by assumption and $\log(1+\epsilon)\prec 1$ by~\ref{L1}, so $\log f \sim \log c >0$.
Suppose $c = 1$, so $\epsilon >0$.
Then~\ref{H1} gives
\[
(\log f)' = f^\dagger \sim \epsilon' < 0,
\]
so $\log f>0$, since $\log f \in \smallo$.
\end{proof}

\begin{lemma}\label{lem:autoL1}
Let $\log\colon K^>\to K$ be an ordered group embedding that satisfies~\ref{L2}, and suppose that $\log(C^>) = C$.
Then $\log$ satisfies~\ref{L1}, so it is a logarithm on $K$.
\end{lemma}
\begin{proof}
Let $a \in 1+\smallo$ and let $c \in C^>$.
Using that $\log(C^>) = C$, take $d_1,d_2 \in C^>$ with $\log d_1 = -c$ and $\log d_2 = c$.
As $\log$ is order-preserving, we have 
\[
d_1<1<d_2 \implies d_1<a<d_2 \implies -c<\log a <c,
\]
so $\log a\in \smallo$, as $c$ was arbitrary.
\end{proof}

The extension $K^\ell$ relates to logarithms as follows:

\begin{lemma}\label{lem:smallexpo}
Let $K$ be a henselian $H$-field and let $\log$ be a logarithm on $K$.
Then $K = K^\ell$ and $\log a = \ln a$ for all $a \in 1+\smallo$.
\end{lemma}
\begin{proof}
We have $\log(1+\smallo)\subseteq \smallo$ by~\ref{L1}, so
\[
(1+\smallo)^\dagger = \der\log(1+\smallo) \subseteq \der\smallo
\]
and $K = K^\ell$.
For $a \in 1+\smallo$, we have $(\log a)' =a^\dagger = (\ln a)'$, so $\log a-\ln a \in C$.
Thus, $\log a = \ln a$, as both belong to $\smallo$.
\end{proof}

\begin{definition}
A \textbf{logarithmic $H$-field} is a henselian $H$-field $K$ with real closed constant field, equipped with a logarithm $\log$ such that the induced logarithm $C^>\to C$ is surjective.
\end{definition}

By Lemma~\ref{lem:autoL1}, the axiom~\ref{L1} is redundant for logarithmic $H$-fields.

\begin{definition}
An \textbf{exponential $H$-field} is a logarithmic $H$-field $K$ with a surjective logarithm.
\end{definition}

Let $K$ be an exponential $H$-field.
Then $\Gamma$ is divisible since each $a \in K^>$ has an $n$-th root $\exp(n\inv\log a)$ for each $n>0$.
As $C$ is real closed and $K$ is henselian, $K$ is real closed by Fact~\ref{fact:realclosed}.
If an element $a \in K$ can be integrated, then it can be exponentially integrated since $b' = (\exp b)^\dagger$ for all $b \in K$.
To summarize:

\begin{remark}\label{rem:realclosed-log}
Every exponential $H$-field is real closed.
Every exponential $H$-field that is closed under integration is Liouville closed.
In particular, every $\upo$-free newtonian exponential $H$-field is Liouville closed, by Fact~\ref{fact:newtonint}.
\end{remark}

Let $K$ and $L$ be logarithmic $H$-fields.
An $H$-field embedding $\iota\colon K\to L$ is said to \textbf{preserve logarithms} if $\iota(\log f) = \log\iota(f)$ for all $f \in K^>$.
A \textbf{logarithmic $H$-field embedding} is an $H$-field embedding that preserves logarithms.
The notions of a logarithmic $H$-field extension and a logarithmic $H$-subfield are defined analogously.
An \textbf{exponential $H$-field extension} is a logarithmic $H$-field extension that is also an exponential $H$-field; exponential $H$-subfields are defined similarly.
For $\phi \in K^>$, the logarithm on $K$ also makes $K^\phi$ a logarithmic $H$-field.

The map $\ln$ gives us a criterion for checking whether a given $H$-field embedding preserves logarithms.

\begin{lemma}\label{lem:logclosedtest}
Let $K$ and $L$ be logarithmic $H$-fields and let $\iota\colon K \to L$ be an $H$-field embedding.
Let $f,g\in K^>$ with $f \sim g$ and suppose $\iota(\log g) = \log\iota (g)$.
Then $\iota (\log f)= \log\iota(f)$.
\end{lemma}
\begin{proof}
Take $\epsilon \in \smallo$ such that $f = g(1+\epsilon)$.
Then
\[
\log f = \log g + \log(1+\epsilon) = \log g + \ln(1+\epsilon)
\]
by Lemma~\ref{lem:smallexpo}.
We need only show that $\iota \ln(1+\epsilon)=\ln\iota(1+\epsilon)$.
To see this, note that $\iota\ln (1+\epsilon)\in \smallo_L$ and that
\[
\iota(\ln (1+\epsilon))' = \iota(1+\epsilon)^\dagger = (\ln\iota(1+\epsilon))'.
\]
Thus, $\iota\ln (1+\epsilon)$ and $\ln\iota(1+\epsilon)$ differ by a constant, so they must be equal as both are in $\smallo_L$.
\end{proof}

The proof of the following proposition borrows from~\cite[Proposition~6.8]{Kap23} and~\cite[Proposition~3.4]{Kap26}.
It will allow us to upgrade a number of $H$-field embeddings from~\cite{ADH17} to logarithmic $H$-field embeddings.

\begin{proposition}\label{prop:logclosedtest2}
Let $K$ be a logarithmic $H$-field and let $M$ be an $H$-field extension of $K$ with real closed constant field.
Suppose we have the following:
\begin{enumerate}[label=\normalfont(\roman*), align = right]
\item A surjective logarithm $C_M^>\to C_M$ extending the induced logarithm on $C$.
\item A family of elements $(a_i)_{i \in I}$ in $M^>$ such that $v(a_i) \neq 0$ for each $i$ and
\[
\Gamma_M = \Gamma \oplus \bigoplus_{i \in I}\Z v a_i.
\]
\item A family of elements $(b_i)_{i \in I}$ in $M$ such that $b_i' = a_i^\dagger$ for each $i$.
\end{enumerate}
Then there is a unique logarithm on $M^\ell$ extending the logarithms on $K$ and $C_M$ such that $\log a_i = b_i$ for each $i \in I$.
With this logarithm, $M^\ell$ is a logarithmic $H$-field extension of $K$.
Moreover, if $L$ is a logarithmic $H$-field extension of $K$ and $\iota\colon M \to L$ is an $H$-field embedding over $K$, then the unique $H$-field embedding $M^\ell \to L$ extending $\iota$ preserves logarithms if and only if $\iota$ commutes with the logarithm on $C_M$ and $\iota(b_i) = \log \iota(a_i)$ for each $i \in I$.
\end{proposition}
\begin{proof}
Since $M^\ell$ is an immediate extension of $M$, the conditions on $M$ also hold for $M^\ell$, so we may assume that $M = M^\ell$.
Let $f \in M^>$.
Our assumption on $\Gamma_M$ gives
\[
f = cg(1+\epsilon)\prod_{i \in I}a_i^{k_i}
\]
for some $c \in C_M^>$, some $g \in K^>$, some $\epsilon \in \smallo_M$, and some family $(k_i)_{i \in I}$ of integers where only finitely many $k_i$ are nonzero.
We set
\[
\log f \coloneqq \log c+\log g + \ln(1+\epsilon) + \sum_{i \in I}k_ib_i,
\]
where $\log c$ is computed using the given logarithm on $C_M$.
It is routine to check that this does not depend on the choice of $c$, $g$, and $\epsilon$.
Then $\log(1+\smallo) = \ln(1+\smallo) \subseteq \smallo$, so $\log$ satisfies~\ref{L1}.
For~\ref{L2}, note that
\[
(\log f)' = g^\dagger + (1+\epsilon)^\dagger + \sum_{i \in I}k_ia_i^\dagger = f^\dagger.
\]
It is straightforward to show that $\log$ is a group morphism, so $\log$ is a logarithm by Lemma~\ref{lem:autoorder}.
As $M = M^\ell$ is henselian, $C_M$ is real closed, and the logarithm on $C_M$ is surjective, $(M,\log)$ is a logarithmic $H$-field.
Uniqueness follows easily from Lemma~\ref{lem:smallexpo}.

Now let $L$ be a logarithmic $H$-field extension of $K$ and let $\iota\colon M \to L$ be an $H$-field embedding over $K$.
We continue to assume that $M = M^\ell$ and we assume that $\iota$ commutes with the logarithm on $C_M$ and that $\iota(b_i) = \log \iota(a_i)$ for each $i \in I$.
We need to show that $\iota(\log f) = \log \iota(f)$ where $f$ is as above.
Using the fact that $f\sim cg\prod_{i \in I}a_i^{k_i}$ and Lemma~\ref{lem:logclosedtest}, we may assume that $f = cg\prod_{i \in I}a_i^{k_i}$.
Then
\[
\log f= \log c+\log g + \sum_{i \in I}k_i\log a_i,
\]
so $\iota(\log f) = \log\iota(f)$ by our assumptions on $\iota$.
The other implication, that if $\iota$ preserves logarithms then it commutes with the logarithm on $C_M$ and satisfies $\iota(b_i) = \log \iota(a_i)$, is clear.
\end{proof}

The conditions on $C_M$ and $\Gamma_M$ in the above proposition are always satisfied when $M$ is an immediate $H$-field extension of $K$:

\begin{corollary}\label{cor:immlogH}
Let $K$ be a logarithmic $H$-field and let $M$ be an immediate $H$-field extension of $K$.
Then there is a unique logarithm $\log$ on $M^\ell$ extending the logarithm on $K$.
With this logarithm, $M^\ell$ is a logarithmic $H$-field extension of $K$ and if $L$ is also a logarithmic $H$-field extension of $K$ and $\iota\colon M \to L$ is an $H$-field embedding over $K$, then the unique $H$-field embedding $M^\ell \to L$ extending $\iota$ preserves logarithms.
\end{corollary}

Let $C^*$ be an ordered field extension of $C$, and consider the field $K(C^*)$, where $K$ and $C^*$ are linearly disjoint over $C$.
By~\cite[Lemma~4.6.21 and Propositions~10.5.15 and~10.5.16]{ADH17}, there are a unique derivation, valuation, and ordering on $K(C^*)$ that make $K(C^*)$ an $H$-field extension of $K$ with constant field $C^*$.
We always view $K(C^*)$ as an $H$-field in this way; as an $H$-field, $K(C^*)$ has value group $\Gamma$.
If $L$ is an $H$-field extension of $K$, then these uniqueness statements show that any ordered field embedding $C^*\to C_L$ over $C$ extends uniquely to an $H$-field embedding $K(C^*)\to L$ over $K$.
Combining these facts with Lemma~\ref{lem:ellext} and Proposition~\ref{prop:logclosedtest2}, we have the following:

\begin{corollary}\label{cor:addconstants}
Let $K$ be a logarithmic $H$-field, let $C^*$ be a real closed field extension of $C$, and let $\log^*$ be a surjective logarithm on $C^*$ extending the logarithm on $C$.
Then there is a unique logarithm on $K(C^*)^\ell$ that extends the logarithms on $K$ and $C^*$.
With this logarithm, $K(C^*)^\ell$ is a logarithmic $H$-field extension of $K$ with value group $\Gamma$ and constant field $C^*$.
If $L$ is a logarithmic $H$-field extension of $K$ and $\iota\colon C^*\to C_L$ is an ordered field embedding over $C$ that commutes with $\log^*$, then $\iota$ extends uniquely to a logarithmic $H$-field embedding $K(C^*)^\ell\to L$ over $K$.
\end{corollary}

\section{Extensions of logarithmic \texorpdfstring{$H$}{H}-fields}\label{sec:exts}
In this section, we prove a variety of extension and embedding results for use in Section~\ref{sec:mc}.
\emph{For the remainder of this section, $K$ is a logarithmic $H$-field with logarithm $\log$.}

\subsection{The Newton-exponential closure}
Here we show that any $\upo$-free logarithmic $H$-field has a minimal $\upo$-free newtonian exponential $H$-field extension.
We begin with the newtonization.

\begin{corollary}\label{cor:newtonization}
Let $K$ be $\upo$-free.
Then there is a unique logarithm on $K^{\nt}$ extending the logarithm on $K$.
If $L$ is a newtonian logarithmic $H$-field extension of $K$, then there is a logarithmic $H$-field embedding $K^{\nt}\to L$ over $K$.
\end{corollary}
\begin{proof}
We have $K^{\nt} = (K^{\nt})^\ell$ by Fact~\ref{fact:newtonint}.
Since $K^{\nt}$ is an immediate extension of $K$, the logarithm on $K$ extends uniquely to a logarithm on $K^{\nt}$ by Corollary~\ref{cor:immlogH}.
The embedding property follows by Fact~\ref{fact:ntexist} and Corollary~\ref{cor:immlogH}.
\end{proof}

\begin{lemma}\label{lem:realclosure}
Set $K^{\rc,\ell}\coloneqq (K^{\rc})^\ell$.
Then $K^{\rc,\ell}$ is a real closed $H$-field with constant field $C$ and there is a unique logarithm on $K^{\rc,\ell}$ extending the logarithm on $K$.
If $L$ is also a real closed logarithmic $H$-field extension of $K$, then there is a unique logarithmic $H$-field embedding $K^{\rc,\ell}\to L$ over $K$.
\end{lemma}
\begin{proof}
Set $M\coloneqq K^{\rc,\ell}$.
Then $M$ is an immediate extension of the $H$-field $K^{\rc}$, so $\Gamma_M = \Q\Gamma$ and $C_M = C^{\rc} = C$ since $C$ is assumed to be real closed.
Thus, $M$ is real closed since $M$ is henselian; see Fact~\ref{fact:realclosed}.
Let $f \in M^>$, and take $c \in C^>$, $g \in K^>$ and $q \in \Q$ with $f \sim cg^q$.
Take also $\epsilon \in \smallo_M$ with $f = cg^q(1+\epsilon)$ and set
\[
\log f \coloneqq \log c + q\log g + \ln(1+\epsilon).
\]
It is routine to check that $\log$ is indeed a logarithm on $M$ extending the logarithm on $K$ and that this extension is unique.
Now let $L$ be a real closed logarithmic $H$-field extending $K$.
Then by the universal property of the real closure there is a unique $H$-field embedding $K^{\rc}\to L$ over $K$ and this in turn extends uniquely to an embedding $\iota\colon M\to L$.
To see that $\iota$ preserves logarithms, let $c \in C^>$, $g \in K^>$, and $q \in \Q$ and note that
\[
\iota \log (cg^q) = \log c+q\log g = \log\iota(cg^q).
\]
We are done in light of Lemma~\ref{lem:logclosedtest}.
\end{proof}

Finally, we deal with adding exponentials.

\begin{lemma}\label{lem:addexp1}
Suppose $K$ is real closed.
Let $a\in K\setminus \log(K^>)$ with $a \prec 1$.
Then $K$ has an immediate $H$-field extension $K(f)$ with $f \sim 1$ and $f^\dagger = a'$.
There is a unique logarithm on $K(f)^\ell$ extending the logarithm on $K$, and this logarithm satisfies $\log f = a$.
If $L$ is also a logarithmic $H$-field extension of $K$ and $a \in \log(L^>)$, then there is a unique logarithmic $H$-field embedding $K(f)^\ell\to L$ over $K$.
\end{lemma}
\begin{proof}
We claim that $a'\not\in(K^\times)^\dagger$.
Suppose towards a contradiction that $a' = b^\dagger$ for $b \in K^\times$.
Then $b \asymp 1$ since $va' \in (\Gamma^>)'$, so by replacing $b$ with $cb$ for some $c \in C^\times$, we may assume that $b\sim 1$.
Lemma~\ref{lem:smallexpo} gives us that $a = \ln b = \log b$, contradicting our assumption that $a\not\in\log(K^>)$.

With this claim out of the way, we may apply~\cite[Lemma~10.5.18]{ADH17} with $a'$ in place of $s$ to get an immediate $H$-field extension $K(f)$ of $K$ where $f \sim 1$ and $f^\dagger = a'$.
By Corollary~\ref{cor:immlogH} there is a unique logarithm $\log$ on $K(f)^\ell$ extending the logarithm on $K$.
We have $\log f = \ln f = a$, since $a' = f^\dagger$.
Now let $L$ be a logarithmic $H$-field extension of $K$ with $a \in \log(L^>)$.
Then $\exp a \in 1+\smallo_L$, and $\exp a$ is transcendental over $K$, as $K$ is real closed.
Since $(\exp a)^\dagger = a'$,~\cite[Lemma~10.5.18]{ADH17} gives us a unique $H$-field embedding $K(f) \to L$ that sends $f$ to $\exp a$.
This in turn extends uniquely to an $H$-field embedding $K(f)^\ell \to L$, which preserves logarithms by Corollary~\ref{cor:immlogH}.
Since any \emph{logarithmic} $H$-field embedding $K(f)^\ell\to L$ must send $f$ to $\exp a$, this embedding is unique as claimed.
\end{proof}

\begin{lemma}\label{lem:addexp2}
Suppose $K$ is real closed.
Let $a \in K$ with $a- \log b \succ 1$ for all $b \in K^>$.
Then $K$ has an $H$-field extension $K(f)$ with constant field $C$ and with $f^\dagger = a'$.
There is a unique logarithm on $K(f)^\ell$ that extends the logarithm on $K$ and satisfies $\log f = a$.
If $L$ is also a logarithmic $H$-field extension of $K$ and $a \in \log(L^>)$, then there is a unique logarithmic $H$-field embedding $K(f)^\ell\to L$ over $K$.
\end{lemma}
\begin{proof}
We may assume that $a <0$.
Then $a'<0$ as well since $a\succ 1$.
We claim that $v(a'-b^\dagger)\in \Psi^{\downarrow}$ for all $b \in K^\times$.
Let $b$ be given; we may assume that $b>0$.
As $v(a-\log b)<0$, we have $v(a'-b^\dagger) = v((a-\log b)') \in (\Gamma^<)'\subseteq \Psi^\downarrow$, as claimed.

In light of this claim, we may apply~\cite[Lemma~10.5.20]{ADH17} with $a'$ in place of $s$ to get an $H$-field $K(f)$ extending $K$ where
\[
f^\dagger = a',\qquad 0<f \prec 1,\qquad\Gamma_{K(f)} = \Gamma \oplus\Z vf,\qquad C_{K(f)} = C.
\]
By Proposition~\ref{prop:logclosedtest2} there is a unique logarithm on $K(f)^\ell$ with $\log f = a$.
Now let $L$ be a logarithmic $H$-field extension of $K$ with $a \in \log(L^>)$.
Note that $\exp a\in L$ is transcendental over $K$, as $K$ is real closed.
As $(\exp a)^\dagger = a'$, there is a unique $H$-field embedding $K(f) \to L$ sending $f$ to $\exp a$ by~\cite[Lemma~10.5.20]{ADH17}.
This in turn extends uniquely to an $H$-field embedding $\iota\colon K(f)^\ell \to L$.
Since $a = \log\iota (f)$, Proposition~\ref{prop:logclosedtest2} gives that $\iota$ preserves logarithms.
Since any logarithmic $H$-field embedding $K(f)^\ell\to L$ must send $f$ to $\exp a$, the map $\iota$ is unique as claimed.
\end{proof}

We now show that each logarithmic $H$-field has a minimal exponential $H$-field extension.

\begin{corollary}\label{cor:expclosure}
The logarithmic $H$-field $K$ has an exponential $H$-field extension $K^e$ with constant field $C$ and with the following property:\ for any exponential $H$-field extension $L$ of $K$, there is a unique logarithmic $H$-field embedding $K^e \to L$ over $K$.
\end{corollary}
\begin{proof}
We define an \emph{$e$-tower} on $K$ to be an increasing chain of logarithmic $H$-fields $(K_\mu)_{\mu\leq\nu}$ where
\begin{enumerate}[label=\normalfont(\alph*)]
\item $K_0= K$;
\item $K_\mu = \bigcup_{\lambda<\mu}K_\lambda$ for each limit ordinal $\mu$;
\item if $K_\mu$ is not real closed, then $K_{\mu+1} = K_\mu^{\rc,\ell}$ equipped with the logarithm in Lemma~\ref{lem:realclosure};
\item if $K_\mu$ is real closed, then $K_{\mu+1} = K_\mu(f_\mu)^\ell$ with $\log f_\mu = a_\mu$, where $a_\mu\in K_\mu\setminus\log(K_\mu^>)$ and either $a_\mu \prec 1$ and $K_\mu(f_\mu)^\ell$ is as in Lemma~\ref{lem:addexp1}, or $a_\mu- \log b\succ 1$ for all $b \in K_\mu^>$ and $K_\mu(f_\mu)^\ell$ is as in Lemma~\ref{lem:addexp2}.
\end{enumerate}
We let $(K_\mu)_{\mu\leq\nu}$ be a maximal $e$-tower on $K$; such a tower exists because each $K_\mu$ is a Liouville extension of $K$, as defined in~\cite[Section~10.6]{ADH17}, so $|K_\mu| = |K|$~\cite[Lemma~10.6.8]{ADH17}.
To see that $K_\nu$ is an exponential $H$-field, we first note that $K_\nu$ is real closed and we let $a \in K_\nu$.
If $a- \log b \succ 1$ for all $b \in K_\nu^>$, then $a \in \log(K_\nu^>)$ by maximality.
Suppose we have $b \in K_\nu^>$ with $a-\log b \preceq 1$.
Since $\log(C^>) = C$, we may take $c \in C^>$ with 
\[
a-\log(cb) = a-\log b-\log c \prec 1.
\]
Again by maximality, $a-\log(cb) \in \log(K_\nu^>)$, so $a \in \log(K_\nu^>)$, as $\log(K_\nu^>)$ is an additive subgroup of $K_\nu$.
Set $K^e\coloneqq K_\nu$.
Then the desired universal property for $K^e$ follows from the universal properties in Lemmas~\ref{lem:realclosure},~\ref{lem:addexp1}, and~\ref{lem:addexp2}.
\end{proof}

Note that $K^e$ is $\d$-algebraic over $K$.
If $K$ is $\upo$-free, then a \textbf{Newton-exponential closure of $K$} is by definition a newtonian exponential $H$-field extension of $K$ that admits a logarithmic $H$-field embedding over $K$ into any other newtonian exponential $H$-field extension of $K$.
Alternating Corollaries~\ref{cor:expclosure} and~\ref{cor:newtonization}, we see that Newton-exponential closures exist.

\begin{proposition}\label{prop:ntexpclosure}
Let $K$ be $\upo$-free.
Then $K$ has a Newton-exponential closure $K^{\nt,e}$ that is $\d$-algebraic over $K$ with constant field $C$.
\end{proposition}

Since any other Newton-exponential closure of $K$ embeds into $K^{\nt,e}$ over $K$, we see that \emph{any} Newton-exponential closure of $K$ is $\d$-algebraic over $K$ (and hence $\upo$-free) with constant field $C$.
By Remark~\ref{rem:realclosed-log}, any Newton-exponential closure is Liouville closed, so by Fact~\ref{fact:altnl}, the underlying $H$-field of any Newton-exponential closure of $K$ is a Newton-Liouville closure of the $H$-field $K$.
Since any two Newton-Liouville closures of $K$ are isomorphic over $K$ by Fact~\ref{fact:nlexist}, we see that any Newton-Liouville closure $L$ of $K$ admits a logarithm that makes $L$ a Newton-exponential closure of $K$.
The following corollary is not used anywhere, but it may be of interest.

\begin{corollary}\label{cor:ntexpclosure}
Let $K$ be $\upo$-free and let $K^{\nt,e}$ be a Newton-exponential closure of $K$.
Then the only newtonian exponential $H$-subfield of $K^{\nt,e}$ containing $K$ is $K^{\nt,e}$ itself.
In particular, any two Newton-exponential closures of $K$ are isomorphic over $K$.
\end{corollary}
\begin{proof}
Let $L$ be a newtonian exponential $H$-subfield of $K^{\nt,e}$ containing $K$.
Then $L$ is $\d$-algebraic over $K$, and so $\upo$-free, newtonian and Liouville closed.
As $K^{\nt,e}$ is a $\d$-algebraic $H$-field extension of $L$ with the same constant field, $K^{\nt,e} = L$ by~\cite[Theorem~16.0.3]{ADH17}.
The uniqueness up to isomorphism follows from the embedding property and the first part of the corollary.
\end{proof}

\subsection{Constructing \texorpdfstring{$\upo$}{ω}-free logarithmic \texorpdfstring{$H$}{H}-field extensions}
In this subsection, we show that any logarithmic $H$-field can be extended to an $\upo$-free logarithmic $H$-field.

\begin{lemma}\label{lem:makegap}
$K$ has a logarithmic $H$-field extension with a gap and with constant field $C$.
\end{lemma}
\begin{proof}
By Lemma~\ref{lem:realclosure}, we may assume that $K$ is real closed.
If $K$ has a gap, then we are done, so as $K$ is ungrounded, we may also assume that $K$ has asymptotic integration.
By~\cite[Corollary~11.4.10 and Lemma~10.5.8]{ADH17}, $K$ has a spherically complete immediate $H$-field extension $L$.
Then $L$ is real closed with asymptotic integration and the logarithm on $K$ extends uniquely to a logarithm on $L^\ell = L$ by Corollary~\ref{cor:immlogH}.
Thus, we may assume that $K$ is spherically complete.
Let $(\upl_\rho)$ be a $\upl$-sequence in $K$, as defined in~\cite[Section~11.5]{ADH17}.
As $K$ is spherically complete, there is $\upl\in K$ with $\upl_\rho\leadsto\upl$.
By Fact~\ref{fact:intclosure} together with spherical completeness, we find $a \in K$ with $a' = -\upl$.
Then $a- \log b \succ 1$ for all $b\in K^>$, since $v\big((a-\log b)'\big) = v(\upl+b^\dagger) \in\Psi^{\downarrow}$ by~\cite[Lemma~11.5.13]{ADH17}.
Thus, we may apply Lemma~\ref{lem:addexp2} to get a logarithmic $H$-field $K(f)^\ell$ extending $K$ with $\log f = a$, so $f^\dagger = a' = -\upl$.
By~\cite[Lemma~11.5.14]{ADH17} and the remark following it, $vf$ is a gap in $K(f)$, and hence also in the immediate extension $K(f)^\ell$ of $K(f)$.
\end{proof}

\begin{lemma}\label{lem:gap2wfree}
Let $s \in K^>$ and suppose that $vs$ is a gap in $K$.
Then $K$ has an $\upo$-free logarithmic $H$-field extension $K(x)^\ell_{\upo}$ with constant field $C$ that contains an element $x\succ 1$ with $x' = s$, and that has the following property:\ if $L$ is also a logarithmic $H$-field extension of $K$ containing $x^*\succ 1$ with $(x^*)' = s$, then there is a unique logarithmic $H$-field embedding $K(x)^\ell_{\upo}\to L$ over $K$ that sends $x$ to $x^*$.
\end{lemma}
\begin{proof}
By~\cite[Corollary~10.5.11]{ADH17}, $K$ has an $H$-field extension $K(x)$ with $x>0$, $x \succ 1$, $x'= s$, and the following properties:
\[
C_{K(x)} = C,\qquad \Gamma_{K(x)} = \Gamma\oplus \Z vx,\qquad \Psi_{K(x)} = \Psi\cup\big\{\psi(vx)\big\},\qquad \Psi < \psi(vx).
\]
In particular, $K(x)$ is grounded with $\max\Psi_{K(x)} = v(x^\dagger)$.
Thus $K(x)$ does not itself admit a logarithm.

Let $K(x)_{\upo}$ be the $H$-field extension of $K(x)$ constructed in~\cite[Section~11.7]{ADH17}, with $x$ as the distinguished element, so $K(x)_{\upo} = K(x_0,x_1,\ldots)$ with $x_0 = x$ and $x_{n+1}' = x_n^\dagger$ for each $n$.
Then $K(x)_{\upo}$ is $\upo$-free with constant field $C$ by~\cite[Corollary~11.7.15]{ADH17}, each $x_n$ is positive, and
\[
\Gamma_{K(x)_{\upo}} = \Gamma\oplus\bigoplus_n \Z vx_n.
\]
The extension $K(x)^\ell_{\upo}$ of $K(x)_{\upo}$ is immediate and $\d$-algebraic, so it is $\upo$-free with constant field $C$.
By Proposition~\ref{prop:logclosedtest2}, there is a unique logarithm on $K(x)^\ell_{\upo}$ that extends the logarithm on $K$ and satisfies $\log x_n = x_{n+1}$ for each $n$.

Now let $L$ and $x^*$ be as in the statement of the lemma, so $x^*$ must be positive by~\ref{H1}.
By~\cite[Corollary~10.5.11]{ADH17}, there is an $H$-field embedding $K(x)\to L$ over $K$ that sends $x$ to $x^*$.
Using the embedding property in~\cite[Corollary~10.5.12]{ADH17}, we extend this embedding inductively to an $H$-field embedding $K(x)_{\upo}\to L$ over $K$ that sends $x_n$ to $\log_nx^*$ for each $n$.
This extends uniquely to a logarithmic $H$-field embedding $K(x)^\ell_{\upo}\to L$ by Proposition~\ref{prop:logclosedtest2}.
Any logarithmic $H$-field embedding $K(x)^\ell_{\upo}\to L$ over $K$ that sends $x$ to $x^*$ must send $x_n$ to $\log_nx^*$ for each $n$, so this embedding is unique.
\end{proof}

In Lemma~\ref{lem:gap2wfree}, we chose the integral of $s$ to be infinite.
Using~\cite[Corollary~10.5.10]{ADH17} in place of~\cite[Corollary~10.5.11]{ADH17} and applying the construction in~\cite[Section~11.7]{ADH17} with $y\inv$ in place of $x$, we have another version where $s$ has an infinitesimal integral.

\begin{lemma}\label{lem:gap2wfree2}
Let $s \in K^>$ and suppose that $vs$ is a gap in $K$.
Then $K$ has an $\upo$-free logarithmic $H$-field extension $K(y)^\ell_{\upo}$ with constant field $C$ that contains an element $y\prec 1$ with $y' = -s$, and that has the following property:\ if $L$ is also a logarithmic $H$-field extension of $K$ containing $y^*\prec 1$ with $(y^*)' = -s$, then there is a unique logarithmic $H$-field embedding $K(y)^\ell_{\upo}\to L$ over $K$ that sends $y$ to $y^*$.
\end{lemma}

The following is immediate from Lemmas~\ref{lem:makegap} and~\ref{lem:gap2wfree}.

\begin{corollary}\label{cor:upofreeext}
$K$ has an $\upo$-free logarithmic $H$-field extension with constant field $C$.
\end{corollary}

Corollary~\ref{cor:upofreeext} and Proposition~\ref{prop:ntexpclosure} give us the following:

\begin{corollary}\label{cor:upontexpclosure}
$K$ has an $\upo$-free newtonian exponential $H$-field extension with constant field $C$.
\end{corollary}

\subsection{The cut \texorpdfstring{$\Gamma^<$}{Γ<}}
In this subsection, we assume that $K$ is $\upo$-free and real closed.
We will show that there is a unique way of realizing the cut $\Gamma^<$.
First, we show how to realize the cut $\Psi^\downarrow$.

\begin{proposition}\label{prop:zgap}
Let $z>0$ be an element in a logarithmic $H$-field extension of $K$ with $\Psi^{\downarrow}<vz<(\Gamma^>)'$.
Then $K_z\coloneqq K\langle z,\log z\rangle^\ell$ is a logarithmic $H$-field with constant field $C$ and gap $vz$ and for any logarithmic $H$-field extension $L$ of $K$ and any $z^* \in L^>$ with $\Psi^{\downarrow}<vz^*<(\Gamma^>)'$, there is a unique logarithmic $H$-field embedding $K_z\to L$ over $K$ that sends $z$ to $z^*$.
\end{proposition}
\begin{proof}
Let $L$ be a logarithmic $H$-field extension of $K$ and let $z^*\in L^>$ with $\Psi^\downarrow <vz^*<(\Gamma^>)'$.
We will build a logarithmic $H$-field embedding $K_z \to L$ in stages, first embedding $K\langle z^\dagger\rangle^{\h}$, then $K\langle z^\dagger\rangle^{\h}(\log z)$, then $K\langle \log z\rangle^\ell$, then finally $K_z$.
Note that there is at most one logarithmic $H$-field embedding $K_z\to L$ over $K$ that sends $z$ to $z^*$, since such an embedding must send $\log z$ to $\log z^*$, so we need only show that such an embedding exists.
We include a diagram below, to help keep track of the various extensions involved.
\[
\begin{tikzcd}[column sep=5em, row sep=1.8em]
K_z = K\langle z,\log z\rangle^\ell \arrow[rr, "z\mapsto z^*", near start] \arrow[d, dash] & & L\\
K\langle \log z \rangle^\ell \arrow[urr, bend right=4] \arrow[d, dash] & & \\
K\langle z^\dagger\rangle^{\h}(\log z) \arrow[uurr, bend right=8, "\log z\mapsto \log z^*" sloped, near start] \arrow[d, dash] & & \\
K\langle z^\dagger\rangle^{\h} \arrow[uuurr, bend right=12, "z^\dagger\mapsto (z^*)^\dagger" sloped, near start] \arrow[d, dash] & & \\
K \arrow[uuuurr, bend right=16] & &
\end{tikzcd}
\]
To start,~\cite[Lemma~13.4.11]{ADH17} gives a valued differential field embedding $K\langle z\rangle \to L$ over $K$ that sends $z$ to $z^*$.
By~\cite[Lemma~13.4.10 and Corollary~9.8.6]{ADH17}, $vz$ is a gap in $K\langle z\rangle$, so $K\langle z^\dagger \rangle$ is an immediate extension of $K$ by~\cite[Lemma~1.4]{Kap26}.
Using also~\cite[Lemma~10.5.8]{ADH17}, we see that the embedding $K\langle z\rangle\to L$ restricts to an $H$-field embedding $K\langle z^\dagger \rangle\to L$ sending $z^\dagger$ to $(z^*)^\dagger$.
We extend this to an embedding $K\langle z^\dagger \rangle^{\h}\to L$ using the universal property of the henselization.

We now turn to adjoining $\log z$.
Note that we can adjoin \emph{some} integral of $z^\dagger$ using Fact~\ref{fact:intclosure}, but this does not necessarily give us the freedom to choose which integral to adjoin or what to send it to.
We will instead show that we are in the setting of Fact~\ref{fact:largeint}, so that we can adjoin $\log z$ and send it to $\log z^*$.
To apply this fact, we need to show the following:
\begin{claim*}
For all $y \in K\langle z^\dagger\rangle^{\h}$, we have $v(z^\dagger-y')<(\Gamma^>)'$.
\end{claim*}
\begin{proof}[Proof of Claim]
We will need to use $\upl$-sequences, as defined in~\cite[Section~11.5]{ADH17}.
Let $(\upl_\rho)$ be a $\upl$-sequence in $K$ and set $\upl\coloneqq -z^\dagger$.
Then $\upl$ and $\upl+z$ are pseudolimits of $(\upl_\rho)$ and $(\upl_\rho)$ is of $\d$-transcendental type over $K$ by~\cite[Corollaries~11.5.7 and~13.6.3]{ADH17}.
By~\cite[Lemma~11.4.7]{ADH17}, both $K\langle\upl\rangle$ and $K\langle\upl+z\rangle$ are immediate extensions of $K$, and there is a valued differential field isomorphism $K\langle\upl\rangle\to K\langle\upl+z\rangle$ over $K$ sending $\upl$ to $\upl+z$.
This extends to an isomorphism $\iota\colon K\langle\upl\rangle^{\h}\to K\langle\upl+z\rangle^{\h}$.

Suppose towards a contradiction that there is $y \in K\langle\upl\rangle^{\h}$ with $v(z^\dagger-y') = v(\upl+y')\in (\Gamma^>)'\cup \{\infty\}$.
Set $w \coloneqq y - \iota y$, so $w \in K\langle \upl,\upl+z\rangle^{\h} = K\langle z\rangle^{\h}$ and 
\[
v(z-w') = v(\upl+z+\iota y'- \upl-y') = v(\iota(\upl+y')-(\upl+y'))\geq v(\upl+y').
\]
As $vz<(\Gamma^>)'$, this gives $w'\sim z$, so $v(w') = vz$.
But $vz$ is a gap in $K\langle z\rangle$, hence in $K\langle z\rangle^{\h}$ as well, a contradiction.
\end{proof}
Having established this claim, Fact~\ref{fact:largeint} gives an immediate $H$-field extension $K\langle z^\dagger\rangle^{\h}(a)$ with $a' = z^\dagger$.
The embedding property in this fact allows us to identify $K\langle z^\dagger\rangle^{\h}(a)$ with $K\langle z^\dagger\rangle^{\h}(\log z)$ by sending $a$ to $\log z$, and it also gives an $H$-field embedding $K\langle z^\dagger\rangle^{\h}(\log z)\to L$ sending $\log z$ to $\log z^*$.
Corollary~\ref{cor:immlogH} gives a unique logarithm on $K\langle \log z\rangle^\ell$ extending the logarithm on $K$, and allows us to extend our embedding $K\langle z^\dagger\rangle^{\h}(\log z)\to L$ to a logarithmic $H$-field embedding $K\langle \log z\rangle^\ell\to L$.

Now we will adjoin $z$.
Note that $K\langle \log z\rangle^\ell$ is real closed, as it is an immediate henselian extension of the real closed field $K$.
For $b \in K\langle \log z\rangle^\ell$ with $b>0$, we have $z\not\asymp b$, as $K\langle \log z\rangle^\ell$ is an immediate extension of $K$ and $vz \not\in \Gamma$.
Thus, $\log z-\log b = \log(z/b) \succ 1$ by Lemma~\ref{lem:logasymp}, and we may apply Lemma~\ref{lem:addexp2} to obtain a logarithmic $H$-field extension $K\langle \log z,f\rangle^\ell$ with $\log f = \log z$.
Using the embedding property of this extension, we may identify this extension with $K_z$ by sending $f$ to $z$, and we get a logarithmic $H$-field embedding $K_z \to L$ sending $z$ to $z^*$.
The proof of Lemma~\ref{lem:addexp2} gives $\Gamma_{K_z} = \Gamma\oplus\Z vz$ and $C_{K_z} = C$, and~\cite[Corollary~9.8.6]{ADH17} gives that $vz$ is a gap in $K_z$.
\end{proof}

Using this proposition, we now show how to realize the cut $\Gamma^<$.

\begin{corollary}\label{cor:endelem}
Let $y>0$ be an element in a logarithmic $H$-field extension of $K$ with $\Gamma^<<vy<0$.
Then $K$ has an $\upo$-free logarithmic $H$-field extension $K_y$ with constant field $C$ containing $y$ such that for any logarithmic $H$-field extension $L$ of $K$ and any $y^* \in L^>$ with $\Gamma^<<vy^*<0$, there is a unique logarithmic $H$-field embedding $K_y\to L$ over $K$ that sends $y$ to $y^*$.
\end{corollary}
\begin{proof}
Let $L$ and $y^*$ be given.
As $y,y^*>C$, the axiom~\ref{H1} gives $y',(y^*)'>0$.
We also have $\Psi \subseteq (\Gamma^<)'<vy',v(y^*)'<(\Gamma^>)'$, as $K$ is ungrounded, so Proposition~\ref{prop:zgap} applies to both $y'$ and $(y^*)'$.
Let $K_{y'} \coloneqq K\langle y',\log(y')\rangle^\ell$ be the logarithmic $H$-field extension of $K$ containing $y'$ as in that proposition, so $vy'$ is a gap in $K_{y'}$ and we have a unique logarithmic $H$-field embedding $K_{y'}\to L$ over $K$ that sends $y'$ to $(y^*)'$.
Using Lemma~\ref{lem:gap2wfree}, we construct an $\upo$-free logarithmic $H$-field $K_y$ extending $K_{y'}$ and containing an infinite integral of $y'$.
Using the embedding property in that lemma, we may assume that this integral is $y$, and we obtain a unique logarithmic $H$-field embedding $K_y\to L$ extending our embedding $K_{y'}\to L$ and sending $y$ to $y^*$.
\end{proof}

\section{Value group extensions}\label{sec:vgext}
Throughout this section, $K$ is an $\upo$-free logarithmic $H$-field, $M$ is a logarithmic $H$-field extension of $K$ with $C_M = C$, and $f$ is an element in $M \setminus K$.
We assume that one of the following two conditions holds:
\begin{enumerate}[label=\normalfont(\Alph*), align = right]
\item\label{caseA} $K$ is a newtonian exponential $H$-field, and $K$ is \emph{maximal} in $M$ in the sense that there is no $y \in M\setminus K$ for which $K\langle y \rangle$ is an immediate extension of $K$;
\item\label{caseB} $\Psi$ is downward closed in $\Gamma$ and $f>K$.
\end{enumerate}
Case~\ref{caseA} is needed for model completeness, whereas case~\ref{caseB} is only used in the proof of local o-minimality (Corollary~\ref{cor:lomin}).
The main result of this section is the following:

\begin{proposition}\label{prop:dtranselem}
$K$ has an $\upo$-free logarithmic $H$-field extension $K_f$ that contains $f$, is an immediate extension of $K\langle f \rangle$, and has the following property:\ if $L$ is also a logarithmic $H$-field extension of $K$ and $f^* \in L$ realizes the same cut as $f$ over $K$, then there is a unique logarithmic $H$-field embedding $K_f\to L$ over $K$ that sends $f$ to $f^*$.
\end{proposition}

This proposition reduces to the more technical Proposition~\ref{prop:approxlog} below.
We first set up what we need to state that proposition, then we explain how it is used to prove Proposition~\ref{prop:dtranselem}.

Let $(\Gamma_f,\psi_f)$ denote the asymptotic couple of $K\langle f\rangle$, with $\Psi_f\coloneqq\psi_f(\Gamma_f^{\neq})$.
We claim that the set
\[
\big\{v(f-a):a \in K\big\} \subseteq \Gamma_f
\]
has a greatest element.
In case~\ref{caseA}, this follows from Fact~\ref{fact:immextmap} and maximality of $K$ in $M$, and in case~\ref{caseB}, this set consists of the single element $vf \not\in \Gamma$.

We define a \textbf{best approximation} to such an $f$ to be an element $b \in K$ such that $v(f-b) = \max\big\{v(f-a):a \in K\big\}$.
Note that $b$ is a best approximation to $f$ if and only if $v(f-b) \not\in\Gamma$ (this uses the fact that $C_M = C$).
In case~\ref{caseB}, any element of $K$ is a best approximation to $f$.

We set $f_0 \coloneqq f$, we let $b_0$ be a best approximation to $f_0$ (taking $b_0 \coloneqq 0$ in case~\ref{caseB}), and we set $f_1\coloneqq (f_0-b_0)^\dagger$.
We claim that $f_1$ satisfies the same conditions as $f$, so $f_1 \not\in K$ in case~\ref{caseA} and $f_1>K$ in case~\ref{caseB}.
In case~\ref{caseA}, this follows from the fact that $K$ is Liouville closed (Remark~\ref{rem:realclosed-log}) and $C_M = C$:\ if $f_1 \in K$ then so is $f_0-b_0$, so $f_0 \in K$.
In case~\ref{caseB}, this uses that $\Psi$ is downward closed:\ if $f_1\not>K$, then there is $g \in K^>$ with $g \succ 1$ and $f_1 \prec g^\dagger$, so $f_0 \preceq g$, contradicting $f_0>K$.
Thus in either case, we can apply the argument above to find a best approximation $b_1$ to $f_1$ in $K$ (again taking $b_1\coloneqq 0$ in case~\ref{caseB}).
Continuing in this way we construct a sequence $(f_n)$ of elements in $M\setminus K$ and a sequence $(b_n)$ of elements in $K$ such that $b_n$ is a best approximation to $f_n$ with $b_n = 0$ in case~\ref{caseB} and such that $f_{n+1} = (f_n-b_n)^\dagger$.

Now for each $n>0$, choose an element $a_n \in K^\times$ such that $a_n^\dagger = b_n$ and such that $a_n$ and $f_{n-1}-b_{n-1}$ have the same sign.
In case~\ref{caseA}, we may do this by Liouville closedness, and in case~\ref{caseB}, we simply choose $a_n \coloneqq 1$ for all $n$.
For each $n$, we set $\fm_n\coloneqq a_{n+1}\inv(f_n-b_n)$, so $\fm_n >0$ and $v(\fm_n) \not\in \Gamma$.
Then
\[
f_{n+1}-b_{n+1}= \left( \frac{f_n-b_n}{a_{n+1}}\right)^\dagger,
\]
so $\fm_n^\dagger = a_{n+2}\fm_{n+1}$ for each $n$.
Note that $K\langle f \rangle = K\langle \fm_0\rangle = K(\fm_0,\fm_1,\ldots)$.
For each $n$, we set
\[
\beta_n\coloneqq v(\fm_n),\qquad \lambda_n \coloneqq v(\log \fm_n),\qquad \alpha_{n+1}\coloneqq v(a_{n+1}).
\]
As $\beta_n\neq 0$, we have $\lambda_n<0$ by Lemma~\ref{lem:logasymp}.
We recall (and establish) a few basic facts about these values:

\begin{lemma}\label{lem:valuegroupgrowth}\
\begin{enumerate}[label=\normalfont(\roman*), align = right]
\item\label{valuegroupgrowthi} $\Gamma_f = \Gamma \oplus \bigoplus_{n}\Z \beta_n$ (internal direct sum).
\item\label{valuegroupgrowthii} $\beta_n^\dagger = \beta_{n+1}+\alpha_{n+2}\not\in\Gamma$ for each $n$; in particular, $[\beta_n]\not\in [\Gamma]$.
\item\label{valuegroupgrowthiii} $\beta_0^\dagger<\beta_1^\dagger<\beta_2^\dagger<\cdots$; in particular, $[\beta_0]>[\beta_1]>[\beta_2]>\cdots$.
\item\label{valuegroupgrowthiv} $[\Gamma_f] = [\Gamma]\cup \{[\beta_n]:n\in\N\}$; in particular, $\Psi_f = \Psi \cup\{\beta_n^\dagger:n\in\N\}$.
\item\label{valuegroupgrowthv} $K\langle f\rangle$ is $\upo$-free, hence has asymptotic integration.
\item\label{valuegroupgrowthvi} $\Gamma^<$ is cofinal in $\Gamma_f^<$, and so $\Psi$ is cofinal in $\Psi_f$.
\item\label{valuegroupgrowthvii} $[\lambda_n] \geq [\beta_{n+1}]$.
\item\label{valuegroupgrowthviii} $\beta_n^\dagger <\lambda_n^\dagger\leq \beta_{n+1}^\dagger$; in particular, $\lambda_m^\dagger\neq \lambda_n^\dagger$ for $m\neq n$.
\end{enumerate}
\end{lemma}
\begin{proof}
In case~\ref{caseA}, parts~\ref{valuegroupgrowthi}--\ref{valuegroupgrowthiv} and~\ref{valuegroupgrowthvi} are~\cite[Lemma~16.1.2]{ADH17} and~\ref{valuegroupgrowthv} is~\cite[Corollary~16.1.3]{ADH17}.
In case~\ref{caseB}, part~\ref{valuegroupgrowthi} and part~\ref{valuegroupgrowthv} follow from~\cite[Lemma~16.6.10]{ADH17}.
Part~\ref{valuegroupgrowthii} follows immediately from~\ref{valuegroupgrowthi}, since $\alpha_n = 0$ for $n>0$, and part~\ref{valuegroupgrowthiii} follows from~\ref{valuegroupgrowthii} and~\cite[Lemma~16.6.9]{ADH17} by induction on $n$.
Part~\ref{valuegroupgrowthiv} follows easily from~\ref{valuegroupgrowthi}--\ref{valuegroupgrowthiii} in this case, and~\ref{valuegroupgrowthvi} follows from~\ref{valuegroupgrowthiv}, since $\Gamma$ is convex in $\Gamma_f$ by~\cite[Lemma~16.6.10]{ADH17}.
Note that the statements of~\cite[Lemmas~16.6.9 and~16.6.10]{ADH17} assume that $K$ is Liouville closed, but their proofs only use that $\Psi$ is downward closed in $\Gamma$.

The proofs of~\ref{valuegroupgrowthvii} and~\ref{valuegroupgrowthviii} are the same in both cases~\ref{caseA} and~\ref{caseB}.
As $\lambda_n' = \beta_n^\dagger$ and $\lambda_n<0$, we have
\[
\lambda_n^\dagger = \beta_n^\dagger-\lambda_n>\beta_n^\dagger,
\]
giving one inequality in~\ref{valuegroupgrowthviii}.
The other inequality in~\ref{valuegroupgrowthviii} follows from~\ref{valuegroupgrowthvii}, so it remains to establish~\ref{valuegroupgrowthvii}.
Suppose towards a contradiction that $[\lambda_n]<[\beta_{n+1}]$.
Then Fact~\ref{fact:6.5.4}\ref{6.5.4i} gives
\[
[\lambda_n^\dagger-\beta_{n+1}^\dagger]<[\lambda_n-\beta_{n+1}] = [\beta_{n+1}].
\]
Using that $\lambda_n^\dagger = \beta_n^\dagger-\lambda_n$, we get $[\beta_n^\dagger - \beta_{n+1}^\dagger] < [\beta_{n+1}]$, but
\[
[\beta_n^\dagger - \beta_{n+1}^\dagger] = [\beta_{n+1}-\beta_{n+2}+\alpha_{n+2}-\alpha_{n+3}] \geq [\beta_{n+1}]
\]
by parts~\ref{valuegroupgrowthii} and~\ref{valuegroupgrowthiii}.
\end{proof}

Note that $\lambda_n\in\Gamma_f$ for each $n$ as $\lambda_n' = \beta_n^\dagger \in \Gamma_f$ and $K\langle f\rangle$ has asymptotic integration by Lemma~\ref{lem:valuegroupgrowth}\ref{valuegroupgrowthv}.
We build a sequence $(M_n)$ of henselian $H$-subfields of $M$ as follows:
\[
M_0\coloneqq K\langle f\rangle^{\h},\qquad M_{n+1}\coloneqq M_n(\log \fm_n)^{\h}.
\]
We can now state our main technical result.

\begin{proposition}\label{prop:approxlog}
For all $n$ and all $y \in M_n$, we have
\[
v(\log\fm_n-y) \leq \lambda_n+o(\lambda_n).
\]
\end{proposition}

The remainder of this section is devoted to proving Proposition~\ref{prop:approxlog}.
First, we show how this proposition is used to prove Proposition~\ref{prop:dtranselem}.

\begin{proof}[Proof of Proposition~\ref{prop:dtranselem}]
We set $M_\infty\coloneqq \bigcup_n M_n$, and we take $K_f\coloneqq M_\infty^\ell$.
We first note that $K_f$ is $\upo$-free:\ this is true for $K\langle f\rangle$ by Lemma~\ref{lem:valuegroupgrowth}\ref{valuegroupgrowthv}, and $K_f$ is $\d$-algebraic over $K\langle f\rangle$ (as $(\log\fm_n)' = \fm_n^\dagger \in K\langle f\rangle$), so it holds for $K_f$ by Fact~\ref{fact:dalgupofree}.

Now we verify the immediacy and universality of $K_f$ by an induction argument.
Let $L$ and $f^*$ be as in the statement of Proposition~\ref{prop:dtranselem}.
By~\cite[Proposition~16.1.5]{ADH17} in case~\ref{caseA} and~\cite[Lemma~16.6.10]{ADH17} and its proof in case~\ref{caseB}, there is a unique $H$-field embedding $\iota\colon K\langle f\rangle\to L$ over $K$ that sends $f$ to $f^*$.
We show by induction on $n$ that:
\begin{enumerate}[label=\normalfont(\alph*)]
\item $M_n$ is an immediate extension of $K\langle f \rangle$.
\item There is a unique $H$-field embedding $\iota_n\colon M_n\to L$ that extends $\iota$ and satisfies $\iota_n(\log\fm_m) = \log \iota(\fm_m)$ for all $m<n$.
\end{enumerate}
When $n = 0$, these conditions follow from immediacy and universality of the henselization.
Suppose that these conditions hold for a given $n\geq 0$.
Proposition~\ref{prop:approxlog} gives
\[
v(\log \fm_n-y) \leq \lambda_n+o(\lambda_n)<0
\]
for all $y \in M_n$, so $v(\fm_n^\dagger- M_n') <(\Gamma_{M_n}^>)'$.
Fact~\ref{fact:largeint} and its universal property give that $M_n(\log\fm_n)$ is an immediate extension of $M_n$ and that there is a unique $H$-field embedding $M_n(\log \fm_n)\to L$ that extends $\iota_n$ and sends $\log \fm_n$ to $\log \iota(\fm_n)$.
Again using immediacy and universality of the henselization, this embedding extends uniquely to the immediate extension $M_{n+1} = M_n(\log \fm_n)^{\h}$.
This completes the induction.

Taking the union of the $\iota_n$'s, we obtain an $H$-field embedding $M_\infty\to L$.
Then $M_\infty$ is an immediate extension of $K\langle f\rangle$, so by Lemma~\ref{lem:valuegroupgrowth}\ref{valuegroupgrowthi} and Proposition~\ref{prop:logclosedtest2} (taking $(\fm_n)$ as the family $(a_i)_{i \in I}$ and $(\log\fm_n)$ for $(b_i)_{i \in I}$), we see that the embedding $M_\infty\to L$ extends to a logarithmic $H$-field embedding $K_f\to L$.
This logarithmic $H$-field embedding is uniquely determined by the condition that $f$ be sent to $f^*$.
\end{proof}

\subsection{A useful conjugation}\label{subsec:conj}
In order to prove Proposition~\ref{prop:approxlog}, we will proceed by induction on $n$, and we will need to compositionally conjugate by a well-chosen element, depending on $n$.

Let $\phi\in K^>$ and let $\derdelta\coloneqq\phi\inv\der$ be the derivation of $K^\phi$; recall that $K^\phi$ and $M^\phi$ are logarithmic $H$-fields with the same logarithms as $K$ and $M$.
Define the sequences $(\tilde a_n)_{n>0}$ and $(\tilde b_n)$ as follows:\ set $\tilde b_0 \coloneqq b_0$, $\tilde b_1 \coloneqq \phi\inv b_1$, $\tilde a_1 \coloneqq a_1$, and for $n\geq 2$, set
\[
\tilde b_n \coloneqq \phi\inv (b_n-\phi^\dagger),\qquad \tilde a_n \coloneqq \phi\inv a_n.
\]
Then $\derdelta(\tilde a_n)/\tilde a_n = \tilde b_n$ for $n \geq 1$.
Define $\tilde f_n$ inductively by
\[
\tilde f_0 \coloneqq f_0,\qquad \tilde f_{n+1} \coloneqq \derdelta(\tilde f_n-\tilde b_n)/(\tilde f_n-\tilde b_n).
\]
We claim that $\tilde f_n - \tilde b_n = \phi\inv(f_n-b_n)$ for each $n>0$.
This is clear for $n = 1$, and assuming that it holds for a given $n \geq 1$, we have
\[
\tilde f_{n+1}- \tilde b_{n+1}= \phi\inv(\phi\inv(f_n-b_n))^{\dagger} - \phi\inv(b_{n+1}-\phi^\dagger)= \phi\inv((f_n-b_n)^\dagger -b_{n+1}) = \phi\inv(f_{n+1}-b_{n+1}).
\]
From this, it follows that $\tilde b_n$ is a best approximation to $\tilde f_n$ and that $\fm_n = \tilde a_{n+1}\inv(\tilde f_n-\tilde b_n)$ for each $n$.
That is to say, working with the compositional conjugates $M^\phi$ and $K^\phi$ instead of $M$ and $K$, the towers $(M_n)$ and the sets $v(M_n-\log\fm_n)$ remain unchanged.
Thus, it suffices to prove Proposition~\ref{prop:approxlog} after compositional conjugation.

\begin{convention}\label{conv:normalize}
For the rest of this section, we fix $n$.
Let $\phi \in K^>$ with $v\phi \in \Psi$ and $v\phi> \lambda_m^\dagger$ for each $m\leq n$ (such $\phi$ exists by Lemma~\ref{lem:valuegroupgrowth}\ref{valuegroupgrowthvi} and~\ref{valuegroupgrowthviii}).
We work in $M^\phi$ and $K^\phi$ in place of $M$ and $K$, and with $\tilde f_m,\tilde a_m,\tilde b_m$ in place of $f_m,a_m,b_m$.
Thus $M$ and $K$ have small derivation, and $\lambda_m^\dagger <0$ for all $m\leq n$.
\end{convention}

For $m\leq n$, set
\[
\Delta_m \coloneqq \{\delta \in \Gamma_f:[\delta]<[\lambda_m]\},
\]
so $\Delta_m$ is a convex subgroup of $\Gamma_f$.
Thus to prove Proposition~\ref{prop:approxlog}, it suffices to show that
\[
v(\log\fm_n-M_n)\subseteq (\lambda_n+\Delta_n)^\downarrow.
\]
For $\delta \in \Delta_m$ we have $[\lambda_m+\delta] = [\lambda_m]$, so
\[
(\lambda_m+\delta)' = \lambda_m+\delta+ (\lambda_m+\delta)^\dagger = \lambda_m+\delta+ \lambda_m^\dagger.
\]
By~\cite[Lemma~9.2.10(iv)]{ADH17}, we have $[\lambda_m^\dagger]<[\lambda_m]$, so $\lambda_m^\dagger \in \Delta_m$.
It follows that
\begin{equation}\label{eq:shift}
(\lambda_m+\Delta_m)' = \lambda_m+\Delta_m<0.
\end{equation}

\subsection{The distance from \texorpdfstring{$K\langle f\rangle$}{K⟨f⟩}}\label{subsec:Kf}
In this subsection, we give a first approximation to Proposition~\ref{prop:approxlog} by showing that $v(\log \fm_n- K\langle f \rangle) \subseteq (\lambda_n+\Delta_n)^\downarrow$.
For each $m$, put $K_m \coloneqq K(\fm_0,\ldots,\fm_m)$, so
\[
K \subsetneq K_0\subsetneq K_1\subsetneq\cdots,\qquad \der(K_m) \subseteq K_{m+1},\qquad K\langle f\rangle = \bigcup_mK_m.
\]
We extend this by setting $K_{-1}\coloneqq K$.
We write $\Gamma_m\coloneqq\Gamma_{K_m}$ with $\Gamma_{-1} = \Gamma$ so, by Lemma~\ref{lem:valuegroupgrowth}\ref{valuegroupgrowthi} and~\cite[Corollary~3.1.11]{ADH17}, we have
\begin{equation}\label{eq:Gammam}
\Gamma_m = \Gamma_{m-1}\oplus\Z\beta_m = \Gamma\oplus\bigoplus_{k\leq m}\Z\beta_k\quad\text{ for all }m\geq 0.
\end{equation}
Fix $m\geq0$.
Then $\fm_m$ is transcendental over $K_{m-1}$ by~\eqref{eq:Gammam} and~\cite[Corollary~3.1.11]{ADH17}, so there is a unique derivation $\theta_m$ of $K_m$ that vanishes on $K_{m-1}$ and satisfies $\theta_m(\fm_m) = \fm_m$.
We set
\[
\der_m\coloneqq\der-\fm_m^\dagger\theta_m.
\]
Then $\der_m(\fm_m) = 0$, so $\der_m$ is a derivation on $K_m$ that coincides with $\der$ on $K_{m-1}$.
Hence for $y\in K_m$,
\[
y' = \der_m(y)+\fm_m^\dagger\,\theta_m(y),
\]
with $\der_m(y),\theta_m(y)\in K_m$.

\begin{lemma}\label{lem:theta}
Let $y \in K_m\setminus K_{m-1}$.
Then $v(\theta_m(y)) \not\in \Gamma_{m-1}$ and there is $z\in K_{m-1}$ with $y-z \asymp \theta_m(y)$.
\end{lemma}
\begin{proof}
Take $P = \sum_k P_k\fm_m^k$ and $Q = \sum_kQ_k\fm_m^k$ in $K_{m-1}[\fm_m]^{\neq}$ with $y = P/Q$.
As $\Gamma_m = \Gamma_{m-1} \oplus \Z\beta_m$, we have
\[
vP = \min_k(vP_k+k\beta_m),
\]
and we let $k_P$ be the unique $k$ at which this minimum is attained; likewise for $k_Q$.
Then
\[
vy = vP_{k_P}- vQ_{k_Q}+(k_P-k_Q)\beta_m.
\]
We have $\theta_m(P) - k_PP = \sum_k(k-k_P)P_k\fm_m^k \prec P$ and $\theta_m(Q) - k_QQ\prec Q$, so
\[
\frac{\theta_m(y)}{y}- k_P+k_Q = \frac{\theta_m(P)}{P}- \frac{\theta_m(Q)}{Q}-k_P+k_Q \prec 1.
\]
First suppose $k_P \neq k_Q$.
Then $v(\theta_m(y)/y) = v(k_P-k_Q) = 0$, so $v(\theta_m(y)) = vy \not\in \Gamma_{m-1}$ and we can take $z = 0$.
Now suppose $k_P = k_Q$, and set $z\coloneqq P_{k_Q}/Q_{k_Q} \in K_{m-1}$.
Then $y - z = R/S$, where $R \coloneqq PQ_{k_Q}-P_{k_Q}Q$ and $S \coloneqq Q_{k_Q}Q$.
Defining $k_R$ and $k_S$ as we did above with $P$ and $Q$, we see that $k_S = k_Q$ and that $k_R \neq k_Q$, since the coefficient of $\fm_m^{k_Q}$ in $R$ is zero.
Appealing to the previous case and using that $\theta_m(z) = 0$, we have
\[
v(\theta_m(y)) = v(\theta_m(y-z)) = v(y-z) \not\in \Gamma_{m-1}.\qedhere
\]
\end{proof}

We will first show that $v(\log\fm_n- y) \in (\lambda_n+\Delta_n)^\downarrow$ for all $y \in K_n$, and then we will show how this carries up to $y \in K_m$ for $m>n$.

\begin{lemma}\label{lem:leveln}
Let $y \in K_n$.
Then $v(\log \fm_n-y) \leq \lambda_n$.\end{lemma}
\begin{proof}
Taking derivatives, it suffices to show that $v(\fm_n^\dagger - y')\leq\beta_n^\dagger$.
We have
\[
\fm_n^\dagger-y' = \fm_n^\dagger (1-\theta_n(y))-\der_n(y)
\]
with $\der_n(y),1-\theta_n(y) \in K_n$.
Since $v(\fm_n^\dagger) = \beta_{n+1}+\alpha_{n+2}\not\in \Gamma_n$ by Lemma~\ref{lem:valuegroupgrowth}\ref{valuegroupgrowthii} and~\eqref{eq:Gammam}, we have
\[
v(\fm_n^\dagger-y') \leq \beta_n^\dagger + v(1-\theta_n(y)).
\]
By Lemma~\ref{lem:theta}, either $\theta_n(y) = 0$ or $v(\theta_n(y)) \not\in \Gamma_{n-1}$; in either case, $1-\theta_n(y)\succeq 1$, so $v(\fm_n^\dagger-y') \leq \beta_n^\dagger$ as desired.
\end{proof}

\begin{lemma}\label{lem:stepup}
Let $m\geq n+1$, let $y \in K_m$, and suppose $v(\log\fm_n-y)> (\lambda_n+\Delta_n)^{\downarrow}$.
Then there is $z \in K_{m-1}$ with $v(\log\fm_n-z) > (\lambda_n+\Delta_n)^{\downarrow}$.
\end{lemma}
\begin{proof}
We may assume $y \not\in K_{m-1}$.
We compute:
\[
\fm_n^\dagger-y' = \fm_n^\dagger-\der_m(y)-\fm_m^\dagger\theta_m(y).
\]
As $\der_m(y)$, $\theta_m(y)$, and $\fm_n^\dagger$ lie in $K_m$ and $v\fm_m^\dagger \not\in \Gamma_m$, we have
\[
\beta_m^\dagger + v(\theta_m(y))\geq v(\fm_n^\dagger-y')>\lambda_n+\Delta_n,
\]
where the last inequality follows from~\eqref{eq:shift} and our assumption on $y$.
We claim that $\beta_m^\dagger \in \Delta_n$.
By Lemma~\ref{lem:valuegroupgrowth}\ref{valuegroupgrowthiii} and~\ref{valuegroupgrowthvii} we have $[\beta_m]\leq[\beta_{n+1}]\leq[\lambda_n]$.
If $[\beta_m] = [\lambda_n]$, then $\beta_m^\dagger = \lambda_n^\dagger \in \Delta_n$.
If $[\beta_m]<[\lambda_n]$, then $[\beta_m^\dagger - \lambda_n^\dagger]<[\beta_m-\lambda_n] = [\lambda_n]$ by Fact~\ref{fact:6.5.4}\ref{6.5.4i}, so $\beta_m^\dagger \in \lambda_n^\dagger + \Delta_n = \Delta_n$.
Thus, $v(\theta_m(y))>\lambda_n+\Delta_n$.
Using Lemma~\ref{lem:theta}, take $z\in K_{m-1}$ with $
y-z\asymp \theta_m(y)$.
Then
\[
v(\log\fm_n-z) \geq \min\{v(\log \fm_n-y),v(y-z)\}>\lambda_n+\Delta_n.\qedhere
\]
\end{proof}

By Lemmas~\ref{lem:leveln} and~\ref{lem:stepup}, we have:

\begin{corollary}\label{cor:stepB}
$v(\log\fm_n-K\langle f\rangle) \subseteq (\lambda_n+ \Delta_n)^\downarrow$.
\end{corollary}

\subsection{The proof of Proposition~\texorpdfstring{\ref{prop:approxlog}}{5.3}}\label{subsec:approxlog}
Now we have shown that $v(\log \fm_n- K\langle f \rangle) \subseteq (\lambda_n+\Delta_n)^\downarrow$, and Proposition~\ref{prop:hensel} handles henselizations, so we only need to deal with the adjunction of $\log\fm_m$ for $m <n$.
The next lemma deals with this step.

\begin{lemma}\label{lem:logstep}
Let $m<n$, and let $E\subseteq M_n$ be an $H$-subfield with $K\langle f\rangle\subseteq E$.
Assume that
\[
v(\log\fm_m-E) \subseteq (\lambda_m+\Delta_m)^\downarrow,\qquad v(\log\fm_n-E) \subseteq (\lambda_n+\Delta_n)^{\downarrow}.
\]
Then $v(\log\fm_n-E(\log\fm_m)) = v(\log\fm_n-E)$.
\end{lemma}
\begin{proof}
Put
\[
S_m\coloneqq v(\log \fm_m-E)\subseteq (\lambda_m+\Delta_m)^\downarrow,\qquad S_n\coloneqq v(\log \fm_n-E) \subseteq (\lambda_n+\Delta_n)^{\downarrow}.
\]
We claim that $S_m$ has no largest element.
Note that $E$ is $\upo$-free, as it is a $\d$-algebraic extension of the $\upo$-free $H$-field $K\langle f\rangle$.
In particular, $E$ has asymptotic integration.
Let $e \in E$, so $v(\log \fm_m-e)<0$ and $v(\log \fm_m-e)' = v(\fm_m^\dagger-e') \in \Gamma_E$.
Asymptotic integration gives $d \in E$ with $d' \asymp \fm_m^\dagger-e'$, so $d \asymp \log \fm_m-e$.
As $C_M = C\subseteq E$, we may multiply $d$ by a constant to arrange that $d \sim \log \fm_m-e$.
Then $v(\log \fm_m-(e+d))>v(\log \fm_m-e)$, which proves the claim.

As $\lambda_m = v(\log \fm_m)\in S_m\subseteq(\lambda_m+\Delta_m)^\downarrow$, we have $\gamma^\dagger= \lambda_m^\dagger$ for all $\gamma \in S_m$ with $\gamma \geq \lambda_m$.
Likewise, $\gamma^\dagger = \lambda_n^\dagger$ for all $\gamma \in S_n$ with $\gamma\geq\lambda_n$.
Lemma~\ref{lem:valuegroupgrowth}\ref{valuegroupgrowthviii} gives $\lambda_m^\dagger\neq \lambda_n^\dagger$, so we may use Proposition~\ref{prop:twoints} to get $v(\log \fm_n-E(\log \fm_m)) = S_n$.
\end{proof}

\begin{proof}[Proof of Proposition~\ref{prop:approxlog}]
The proof proceeds by induction on $n$, so let $n$ be given, assume that Convention~\ref{conv:normalize} is in force for $n$, and suppose that
\[
v(\log \fm_m-M_m) \subseteq (\lambda_m+ \Delta_m)^\downarrow
\]
for each $m<n$.
Corollary~\ref{cor:stepB} gives
\[
v(\log\fm_n-y) \in (\lambda_n+ \Delta_n)^\downarrow
\]
for all $y \in K\langle f\rangle$.
By Proposition~\ref{prop:hensel}, the same holds for all $y \in M_0$, since $(\log \fm_n)' = \fm_n^\dagger \in K\langle f\rangle$.
Let $m< n$ and assume that
\[
v(\log\fm_n- y)\in (\lambda_n+ \Delta_n)^\downarrow
\]
for all $y \in M_m$.
Then this holds for all $y \in M_m(\log \fm_m)$ by Lemma~\ref{lem:logstep}, and so also for all $y \in M_{m+1}$ by Proposition~\ref{prop:hensel}.
Thus, this holds for all $y \in M_n$ by induction.
\end{proof}

\section{A model complete theory of logarithmic \texorpdfstring{$H$}{H}-fields}\label{sec:mc}
In this section, we axiomatize a model complete theory of logarithmic $H$-fields.
We characterize the completions of this theory and show that models of this theory are locally o-minimal with stably embedded exponential constant field.
We will also show how our results can be extended to include restricted sine and cosine functions.
Specific models of interest, including $\T_{\exp}$, are discussed in Section~\ref{sec:examples}.

The proofs in this section are very similar to those in~\cite[Section~5]{Kap26}, which in turn adapt the proofs in~\cite[Chapter~16]{ADH17}.
There are some slight modifications here and in~\cite{Kap26} that allow us to characterize completions and describe the induced structure on $C$ without first proving quantifier elimination.
The main ingredient is the following proposition:\ a direct analog of~\cite[Proposition~16.2.3]{ADH17}.

\begin{proposition}\label{prop:bigembedding}
Let $K$ and $L$ be $\upo$-free newtonian exponential $H$-fields and assume the underlying ordered set of $L$ is $|K|^+$-saturated and the cofinality of $\Gamma_L^<$ is greater than $|\Gamma|$.
Let $E$ be an $\upo$-free logarithmic $H$-subfield of $K$ with $C_E = C$ and let $\iota\colon E \to L$ be a logarithmic $H$-field embedding.
Then $\iota$ extends to a logarithmic $H$-field embedding $K\to L$.
\end{proposition}
\begin{proof}
We can assume $E\neq K$, in which case it suffices to show that $\iota$ can be extended to an embedding of an $\upo$-free logarithmic $H$-subfield of $K$ that properly contains $E$.
Using Proposition~\ref{prop:ntexpclosure}, we identify $E^{\nt,e}$ with a logarithmic $H$-subfield of $K$ and extend $\iota$ to a logarithmic $H$-field embedding $E^{\nt,e} \to L$.
Thus, we may assume for the remainder of the proof that $E$ is an $\upo$-free newtonian exponential $H$-subfield of $K$.
In particular, $\Gamma_E$ is divisible.

Suppose $\Gamma_E^<$ is not cofinal in $\Gamma^<$, so there is $y \in K^>$ with $\Gamma_E^<<vy<0$.
By our cofinality assumption on $\Gamma_L^<$ there is $y^* \in L^>$ with $\Gamma_{\iota(E)}^<<vy^*<0$.
Using Corollary~\ref{cor:endelem}, we extend $\iota$ to a logarithmic $H$-field embedding of an $\upo$-free logarithmic $H$-subfield $E_y\subseteq K$ containing $y$.
Thus, we may assume for the remainder of the proof that $\Gamma_E^<$ is cofinal in $\Gamma^<$, so every differential subfield of $K$ containing $E$ is an $\upo$-free $H$-subfield of $K$ by Fact~\ref{fact:dalgupofree}.

Next, suppose $E\langle y\rangle$ is an immediate extension of $E$ for some $y \in K\setminus E$.
Then $E\langle y\rangle^\ell$ is a logarithmic $H$-subfield of $K$ by Corollary~\ref{cor:immlogH}.
The saturation assumption on $L$ gives $y^* \in L$ that realizes the $\iota$-image of the cut over $\iota(E)$ that $y$ realizes over $E$.
By Fact~\ref{fact:immextmap} we may extend $\iota$ to an $H$-field embedding of $E\langle y\rangle$ into $L$ that sends $y$ to $y^*$.
This in turn extends to an $H$-field embedding $E\langle y\rangle^\ell \to L$.
By Corollary~\ref{cor:immlogH} this $H$-field embedding preserves logarithms.

Finally assume that there is no $y\in K\setminus E$ such that $E\langle y\rangle$ is an immediate extension of $E$.
Then we are in case~\ref{caseA} of Section~\ref{sec:vgext}.
Take $f \in K\setminus E$ and, using the saturation assumption on $L$, take $f^*\in L$ that realizes the $\iota$-image of the cut over $\iota(E)$ that $f$ realizes over $E$.
By Proposition~\ref{prop:dtranselem}, we may extend $\iota$ to a logarithmic $H$-field embedding of an $\upo$-free logarithmic $H$-subfield of $K$ containing $f$.
\end{proof}

Here is a useful consequence of Proposition~\ref{prop:bigembedding}.

\begin{corollary}\label{cor:bigembedding}
Let $K$ and $L$ be as in the statement of Proposition~\ref{prop:bigembedding}.
If both $K$ and $L$ have small derivation, then any ordered exponential field embedding $C\to C_L$ extends to a logarithmic $H$-field embedding $K\to L$.
The same is true if both $K$ and $L$ have large derivation.
\end{corollary}
\begin{proof}
We first consider the case that $K$ and $L$ have small derivation.
Since $K$ is $\upo$-free and newtonian, there is $x \in K$ with $x' = 1$ by Fact~\ref{fact:newtonint}.
Additionally, $x \succ 1$ since $K$ has small derivation.
We view $C$ as a logarithmic $H$-subfield of $K$ with trivial derivation and gap 0 and we let $C(x)^\ell_{\upo}$ be the $\upo$-free logarithmic $H$-field extension of $C$ containing $x$ constructed in Lemma~\ref{lem:gap2wfree}, identified with a logarithmic $H$-subfield of $K$ via the embedding property in that lemma.
Since $L$ is also $\upo$-free and newtonian with small derivation, there is $x^* \in L$ with $x^*\succ 1$ and $(x^*)' =1$, so our embedding $C\to C_L$ extends to a logarithmic $H$-field embedding $C(x)^\ell_{\upo}\to L$ that sends $x$ to $x^*$ by Lemma~\ref{lem:gap2wfree}.
This further extends to a logarithmic $H$-field embedding $K\to L$ by Proposition~\ref{prop:bigembedding}.

Now suppose $K$ and $L$ have large derivation.
This time, we find $y \in K$ with $y' = -1$.
Then $y \preceq 1$, since $K$ has large derivation, so by subtracting a constant from $y$, we may assume that $y\prec 1$.
We also find $y^* \in L$ with $y^*\prec 1$ and $(y^*)' = -1$.
This time using Lemma~\ref{lem:gap2wfree2}, we extend our embedding $C\to C_L$ to a logarithmic $H$-field embedding of $C(y)^\ell_{\upo}\subseteq K$ into $L$.
This further extends to a logarithmic $H$-field embedding $K\to L$, again by Proposition~\ref{prop:bigembedding}.
\end{proof}

\subsection{Model completeness and completeness}\label{subsec:mc}
To get a model completeness result, we need to remove the assumption that $C_E = C$ in Proposition~\ref{prop:bigembedding}.
In order to do this, we impose some additional requirements on the constant fields.

A logarithmic $H$-field $K$ is said to have \textbf{real exponential constant field} if its constant field $C$ equipped with $\exp|_{C}$ models $\Th(\R_{\exp})$.
Let $\cL_{\log}\coloneqq \{+,-,\times,0,1,\leq,\preceq,\der,\log\}$, where $\leq$ and $\preceq$ are binary relation symbols and where $\der$ and $\log$ are unary function symbols.
We view each logarithmic $H$-field $K$ as an $\cL_{\log}$-structure in the obvious way, where $\log$ is defined to be identically zero on $K^{\leq}$.
We consider the following $\cL_{\log}$-theories:
\begin{enumerate}[label=\normalfont(\arabic*), align = right]
\item The theory $T_{\log}$ of logarithmic $H$-fields with real exponential constant field.
\item The theory $T^{\upo}_{\log}$ of $\upo$-free logarithmic $H$-fields with real exponential constant field.
\item The theory $T^{\nt}_{\exp}$ of $\upo$-free newtonian exponential $H$-fields with real exponential constant field.
\end{enumerate}

\begin{theorem}\label{thm:mainthm}
The $\cL_{\log}$-theory $T^{\nt}_{\exp}$ is model complete.
It is the model companion of $T_{\log}$ and the model completion of $T^{\upo}_{\log}$.
\end{theorem}
\begin{proof}
By a standard test for model completions~\cite[p.~viii]{Kap21}, we need to show the following:
\begin{enumerate}[label=\normalfont(\roman*), align=right]
\item\label{mainthmi} Every model of $T_{\log}$ extends to a model of $T^{\upo}_{\log}$.
\item\label{mainthmii} Every model of $T^{\upo}_{\log}$ extends to a model of $T^{\nt}_{\exp}$.
\item\label{mainthmiii} Given $K,L \models T^{\nt}_{\exp}$ with a common substructure $E\models T^{\upo}_{\log}$, if $L$ is $|K|^+$-saturated then the inclusion $E\subseteq L$ extends to an embedding $K\to L$.
\end{enumerate}
Then~\ref{mainthmi} is Corollary~\ref{cor:upofreeext} and~\ref{mainthmii} is Proposition~\ref{prop:ntexpclosure}, so we need only show~\ref{mainthmiii}.
Let $K,L,E$ be as in~\ref{mainthmiii}, so $C,C_L,C_E \models\Th(\R_{\exp})$ and $C_L$ is $|C|^+$-saturated.
The theory of $\R_{\exp}$ is model complete by Wilkie's theorem~\cite{Wil96}, so the inclusion $C_E \subseteq C_L$ extends to an ordered exponential field embedding $C\to C_L$.
By Corollary~\ref{cor:addconstants}, this embedding extends to a unique logarithmic $H$-field embedding $E(C)^\ell\to L$ that is the identity on $E$.
Since $E(C)^\ell$ is $\d$-algebraic over $E$, it is $\upo$-free by Fact~\ref{fact:dalgupofree}.
By Proposition~\ref{prop:bigembedding} with $E(C)^\ell$ in place of $E$, this embedding extends to a logarithmic $H$-field embedding $K\to L$.
\end{proof}

We can use Corollary~\ref{cor:bigembedding} and Theorem~\ref{thm:mainthm} to characterize the completions of $T^{\nt}_{\exp}$.
Let $T^{\nt}_{\exp,\dsmall}$ be the theory of $\upo$-free newtonian exponential $H$-fields with real exponential constant field and small derivation.
We define $T^{\nt}_{\exp,\dlarge}$ analogously, but with large derivation instead.
We can go back and forth between these two theories via compositional conjugation:\ for $K \models T^{\nt}_{\exp}$ and $\phi \in K^>$, we have
\[
K^{\phi} \models T^{\nt}_{\exp,\dsmall}\iff v\phi \in (\Gamma^<)'.
\]

\begin{theorem}\label{thm:completions}
The theories $T^{\nt}_{\exp,\dsmall}$ and $T^{\nt}_{\exp,\dlarge}$ are the two completions of $T^{\nt}_{\exp}$.
These two completions, as well as $T^{\nt}_{\exp}$ itself, are decidable relative to $\Th(\R_{\exp})$.
\end{theorem}
\begin{proof}
Let $K,L \models T^{\nt}_{\exp,\dsmall}$ and assume $L$ is $|K|^+$-saturated.
Then $C_L$ is $|C|^+$-saturated, so there is an ordered exponential field embedding $C\to C_L$ since $C$ and $C_L$ are elementarily equivalent; see~\cite[Corollary~B.9.5]{ADH17}.
This in turn extends to an embedding $K\to L$ by Corollary~\ref{cor:bigembedding}.
As $T^{\nt}_{\exp}$ is model complete, this embedding is elementary, so $K$ and $L$ are elementarily equivalent.
This shows completeness of $T^{\nt}_{\exp,\dsmall}$, and the same proof shows completeness of $T^{\nt}_{\exp,\dlarge}$.
For decidability, we need only note that our axiomatizations of these theories are effective, except possibly the ``real exponential constant field'' axiom; see~\cite[Corollary~16.6.3]{ADH17}.
\end{proof}

For the time being, an effective axiomatization of $\Th(\R_{\exp})$ seems to be out of reach; see~\cite{MW96}.
Let us remark on a reduct of $T^{\nt}_{\exp}$ studied in~\cite{Kap26}:

\begin{remark}\label{rem:pow}
Let $K$ be an exponential $H$-field.
Then by~\cite[Proposition~3.2]{Kap26}, the map
\[
(f,c)\mapsto \exp(c\log f)\colon K^>\times C\to K^>
\]
is a constant power map on $K$ in the sense of~\cite{Kap26}.
If $K\models T^{\nt}_{\exp}$, then $K$ is Liouville closed by Remark~\ref{rem:realclosed-log}, so $K$, equipped with this constant power map and with $\exp(1)$ as distinguished constant, is a model of the theory $T^{\nl}_{\pow}$ from~\cite{Kap26}.

Note that Theorem~\ref{thm:mainthm} does not imply the model completeness of $T^{\nl}_{\pow}$, shown in~\cite[Theorem~5.3]{Kap26}.
However, Corollaries~\ref{cor:lomin} and~\ref{cor:stablyembedded} do generalize the corresponding results~\cite[Corollaries~5.10 and~5.11]{Kap26}.
\end{remark}

We end this subsection with an alternative axiomatization of $T^{\nt}_{\exp}$.
Let $K$ be an $H$-field.
We say that $K$ has the \textbf{differential intermediate value property} (DIVP) if for all $P\in K\{Y\}$ and all $f<g \in K$ with
\[
P(f) < 0 < P(g),
\]
there is $y \in K$ with $f<y<g$ and $P(y) = 0$.
In~\cite{ADH22B}, an alternative axiomatization of $T^{\nl}$ was given using the fact that the field $\T_{\operatorname{g}}$ of grid-based transseries has DIVP~\cite{vdH02}.

\begin{fact}[{\cite[Theorem~2.7]{ADH22B}}]\label{fact:DIVP}
A Liouville closed $H$-field is $\upo$-free and newtonian if and only if it has DIVP.
\end{fact}

We now leverage this to give another axiomatization of $T^{\nt}_{\exp}$.

\begin{corollary}\label{cor:DIVP}
An exponential $H$-field with nontrivial derivation is $\upo$-free and newtonian if and only if it has DIVP.
Thus $T^{\nt}_{\exp}$ is axiomatized by the $\cL_{\log}$-theory of exponential $H$-fields with nontrivial derivation, real exponential constant field, and DIVP.
\end{corollary}
\begin{proof}
Let $K$ be an exponential $H$-field with nontrivial derivation.
If $K$ is $\upo$-free and newtonian, then $K$ is Liouville closed by Remark~\ref{rem:realclosed-log}, so $K$ has DIVP by Fact~\ref{fact:DIVP}.
Conversely, suppose $K$ has DIVP.
As the derivation of $K$ is nontrivial, we have $\Gamma\neq\{0\}$, so $K' = K$ by~\cite[Lemma~2.1]{ADH22B}.
Thus, $K$ is Liouville closed by Remark~\ref{rem:realclosed-log} and hence $\upo$-free and newtonian by Fact~\ref{fact:DIVP}.
\end{proof}

\subsection{Some properties of definable sets}
Let $K\models T^{\nt}_{\exp}$, let $K^*$ be a $|K|^+$-saturated elementary extension of $K$, and set $C^*\coloneqq C_{K^*}$.
In this subsection, we show that $K$ is \emph{locally o-minimal} in the sense of~\cite{TV09} and that $C$ is stably embedded as a model of $\Th(\R_{\exp})$.

\begin{corollary}\label{cor:lomin}
Let $X \subseteq K$ be $\cL_{\log}(K)$-definable.
\begin{enumerate}[label=\normalfont(\arabic*), align = right]
\item\label{lomin}For every $b \in K\cup \{+\infty\}$, there is $a \in K$ with $a<b$ such that $(a,b)$ is either contained in or disjoint from $X$.
\item\label{lominC} There is $f \in K$ with $f>C$ such that $\{y \in K:C<y<f\}$ is either contained in or disjoint from $X$.
\end{enumerate}
\end{corollary}
\begin{proof}
We begin by verifying~\ref{lomin} in the case $b = +\infty$; the general case follows by taking fractional linear transformations.
Let $a_1,a_2 \in K^*$ with $a_1,a_2>K$; we need to show that $a_1$ and $a_2$ have the same $\cL_{\log}$-type over $K$.
As $K$ is Liouville closed, $\Psi$ is downward closed in $\Gamma$.
By the remarks before Corollary~\ref{cor:addconstants}, we have $\Gamma_{K(C^*)^\ell} = \Gamma$, so $a_1,a_2>K(C^*)^\ell$ and $\Psi_{K(C^*)^\ell}$ is downward closed in $\Gamma_{K(C^*)^\ell}$.
Moreover, $K(C^*)^\ell$ is $\upo$-free by Fact~\ref{fact:dalgupofree}, so we are in case~\ref{caseB} of Section~\ref{sec:vgext} with $K(C^*)^\ell$ in place of $K$.
Proposition~\ref{prop:dtranselem} yields an isomorphism of $\upo$-free logarithmic $H$-subfields of $K^*$ over $K(C^*)^\ell$ that sends $a_1$ to $a_2$.
As $T^{\nt}_{\exp}$ is the model completion of $T^{\upo}_{\log}$, this is a partial elementary map $K^*\to K^*$, so $a_1$ and $a_2$ have the same $\cL_{\log}$-type over $K$.

Now we turn to~\ref{lominC}.
Let $f_1,f_2 \in K^*$ with $C^*<f_1,f_2<K^{>C}$; again we need to show that $f_1$ and $f_2$ have the same $\cL_{\log}$-type over $K$.
As $\Gamma$ is divisible, we have $\Gamma^<<vf_1,vf_2<0$.
Corollary~\ref{cor:endelem} yields an isomorphism over $K$ of $\upo$-free logarithmic $H$-subfields of $K^*$ that sends $f_1$ to $f_2$.
As before, this is a partial elementary map $K^*\to K^*$, so $f_1$ and $f_2$ have the same $\cL_{\log}$-type over $K$.
\end{proof}

These properties were established with similar proofs for $T^{\nl}$ in~\cite[Section~16.6]{ADH17}.
The proof of the next property is similar to the proof of~\cite[Corollary~5.11]{Kap26}.

\begin{corollary}\label{cor:stablyembedded}
The constant field $C$ of $K$ is stably embedded as a model of $\Th(\R_{\exp})$:\ every subset of $C^n$ that is $\cL_{\log}(K)$-definable is definable with parameters in the ordered exponential field $C$.
\end{corollary}
\begin{proof}
Suppose $c_1,c_2\in (C^*)^n$ have the same type over $C$ in the ordered exponential field $C^*$.
As in the previous corollary, we need to show that $c_1$ and $c_2$ have the same $\cL_{\log}$-type over $K$.
For $i = 1,2$, let $C_i$ be the definable closure of $C$ and the $c_i$ in the ordered exponential field $C^*$.
Since $\Th(\R_{\exp})$ has definable Skolem functions, $C_i$ is an elementary substructure of $C^*$.
As $c_1$ and $c_2$ have the same type over $C$, there is an ordered exponential field isomorphism $\iota\colon C_1\to C_2$ over $C$ that sends $c_1$ to $c_2$.
By Corollary~\ref{cor:addconstants}, $\iota$ extends to a logarithmic $H$-field embedding $K(C_1)^\ell\to K^*$ over $K$.
As $K(C_1)^\ell$ is $\d$-algebraic over $K$, it is $\upo$-free by Fact~\ref{fact:dalgupofree}, so $K(C_1)^\ell\models T^{\upo}_{\log}$ and this embedding is partial elementary.
Thus $c_1$ and $c_2$ have the same $\cL_{\log}$-type over $K$.
\end{proof}

\subsection{Restricted elementary functions}\label{subsec:trig}
In this subsection, we examine $H$-fields with restricted sine and cosine functions.
Let $K\models T^{\nt}_{\exp}$ and let $[-1,1]_K$ denote the closed interval $\{x\in K:-1\leq x \leq 1\}$.
\textbf{Restricted trigonometric functions on $K$} are functions $\sin\colon [-1,1]_K\to K$ and $\cos\colon [-1,1]_K\to K$ such that
\begin{enumerate}[label=\normalfont(RT\arabic*), align = right]
\item\label{RT1} $\sin (a+b) = \sin a \cos b+\cos a \sin b$ and $\cos( a+b) = \cos a \cos b - \sin a \sin b$ for all $a,b \in [-1,1]_K$ with $a+b \in [-1,1]_K$;
\item\label{RT2} $(\sin a)' = (\cos a)a'$ and $(\cos a)' = -(\sin a)a'$ for all $a \in [-1,1]_K$;
\item\label{RT3} $\sin(\smallo)\subseteq \smallo$ and $\cos(\smallo)\subseteq 1+\smallo$.
\end{enumerate}

Let $\sin$ and $\cos$ be restricted trigonometric functions on $K$.
Then $\sin$ and $\cos$ remain restricted trigonometric functions in any compositional conjugate of $K$.
For each $c \in [-1,1]_C$ we have $(\sin c)' = (\cos c)' = 0$, so $\sin c,\cos c \in C$.
With~\ref{RT3}, this gives us $\sin(0) = 0$ and $\cos(0) = 1$.
The restrictions of $\sin$ and $\cos$ to $\smallo$ are definable in the underlying $H$-field of $K$:\ let $a \in \smallo^{\neq}$ and consider the linear differential equation
\begin{equation}\label{eq:sincos}
-(a')^3Y+a''Y'-a'Y'' = 0.
\end{equation}
Then $\sin a$ is the unique zero of~\eqref{eq:sincos} in $a(1+\smallo)$ and $\cos a$ is the unique zero of~\eqref{eq:sincos} in $1+a\smallo$~\cite[Lemma~8.8]{Kap21}; see also~\cite[Lemma~1.2.17]{ADH24B}.

Let $\R_{\elem}$ be the expansion of $\R_{\exp}$ by restricted sine and cosine functions, so $\Th(\R_{\elem})$ is model complete and o-minimal by van den Dries and Miller~\cite[Corollary~6.12]{DM94}.
Here the subscript stands for ``elementary functions''.
Let $\cL_{\elem}\coloneqq \cL_{\log}\cup\{\sin,\cos\}$, where $\sin$ and $\cos$ are unary function symbols, interpreted as identically zero off $[-1,1]_K$.
Let $T^{\nt}_{\elem}$ be the $\cL_{\elem}$-theory of $\upo$-free newtonian exponential $H$-fields $K$ equipped with restricted trigonometric functions $\sin$ and $\cos$, such that the expansion of $C$ by $\exp|_C$, $\sin|_{[-1,1]_C}$ and $\cos|_{[-1,1]_C}$ is a model of $\Th(\R_{\elem})$.

\begin{lemma}\label{lem:trigembeddingprep}
Let $K,L \models T^{\nt}_{\elem}$, let $\iota\colon K\to L$ be a logarithmic $H$-field embedding, and suppose that
\[
\iota(\sin c) = \sin \iota(c),\qquad \iota(\cos c) = \cos \iota(c)
\]
for all $c \in [-1,1]_C$.
Then $\iota$ is an $\cL_{\elem}$-embedding.
\end{lemma}
\begin{proof}
For $f\in[-1,1]_K\setminus C$, write $f = c+a$ with $c\in[-1,1]_C$ and $a\in\smallo^{\neq}$; by~\ref{RT1} and our assumption on $\iota$, it suffices to show that $\iota(\sin a) = \sin\iota(a)$ and $\iota(\cos a) = \cos\iota(a)$.
This follows from definability of $\sin$ and $\cos$ on the maximal ideal $\smallo$, exactly as in the proof of~\cite[Corollary~8.11]{Kap21}.
\end{proof}

\begin{theorem}\label{thm:trigmainthm}\
\begin{enumerate}[label=\normalfont(\arabic*), align = right]
\item\label{trigmainthmi} $T^{\nt}_{\elem}$ is model complete.
\item\label{trigmainthmii} $T^{\nt}_{\elem}$ has two completions:\ one with small derivation and one with large derivation.
\item\label{trigmainthmiii} Let $K\models T^{\nt}_{\exp}$ and suppose that $C$ admits an expansion to a model of $\Th(\R_{\elem})$.
Then $K$ admits a unique expansion to a model of $T^{\nt}_{\elem}$ extending this expansion of $C$.
\item\label{trigmainthmiv} Let $K \models T^{\nt}_{\elem}$ and let $X \subseteq K$ be $\cL_{\elem}(K)$-definable.
\begin{enumerate}[label=\normalfont(4\alph*), align = right]
\item\label{trigmainthmiva} For $b \in K\cup \{+\infty\}$, there is $a\in K$ with $a<b$ such that $(a,b)$ is either contained in or disjoint from $X$.
\item\label{trigmainthmivb} There is $g \in K$ with $g>C$ such that $\{y \in K:C<y<g\}$ is either contained in or disjoint from $X$.
\end{enumerate}
\item\label{trigmainthmv} Let $K \models T^{\nt}_{\elem}$.
Then the constant field $C$ of $K$ is stably embedded as a model of $\Th(\R_{\elem})$.
\end{enumerate}
\end{theorem}
\begin{proof}
For~\ref{trigmainthmi}, let $K$, $L$, and $E$ be models of $T^{\nt}_{\elem}$ where $E$ is a common $\cL_{\elem}$-substructure of $K$ and $L$ and where $L$ is $|K|^+$-saturated.
By the model completeness of $\Th(\R_{\elem})$, we extend the inclusion $C_E \subseteq C_L$ to an embedding $C\to C_L$ that respects the exponential and restricted trigonometric functions on $C$.
As in the proof of Theorem~\ref{thm:mainthm}, this extends to a logarithmic $H$-field embedding $K\to L$, which is even an $\cL_{\elem}$-embedding in light of Lemma~\ref{lem:trigembeddingprep}.
Now~\ref{trigmainthmi} follows from~\cite[Corollary~B.10.4]{ADH17}.

The proof of~\ref{trigmainthmii} is the same as the proof of Theorem~\ref{thm:completions}:\ let $K,L \models T^{\nt}_{\elem}$ have small derivation, assume that $L$ is $|K|^+$-saturated, take an embedding $C\to C_L$ that respects the exponential and restricted trigonometric functions on $C$, and extend this to a logarithmic $H$-field embedding $K\to L$ using Corollary~\ref{cor:bigembedding}.
This embedding is an $\cL_{\elem}$-embedding by Lemma~\ref{lem:trigembeddingprep}, so it is elementary by~\ref{trigmainthmi}.
The case of large derivation is the same.

One can adapt the proof of~\cite[Proposition~8.13]{Kap21} to establish~\ref{trigmainthmiii}:\ first note that for any $a \in \smallo$, there are solutions to~\eqref{eq:sincos} in $a(1+\smallo)$ and $1+a\smallo$ (this can be deduced from newtonianity explicitly, or it can be established by transfer, as shown in the proof of~\cite[Proposition~8.13]{Kap21}).
These solutions are unique, so we may use the formula in~\ref{RT1} to extend the restricted trigonometric functions on $C$ to functions on all of $[-1,1]_K$.

We prove~\ref{trigmainthmiva} using the proof of Corollary~\ref{cor:lomin}\ref{lomin}; the proof of~\ref{trigmainthmivb} is similar.
It suffices to show that there is a unique $\cL_{\elem}$-type over $K$ at $+\infty$.
Fix a $|K|^+$-saturated elementary extension $K^*$ of $K$ and $a_1,a_2 \in K^*$ with $a_1,a_2>K$.
The proof of Corollary~\ref{cor:lomin}\ref{lomin} gives an isomorphism of $\upo$-free logarithmic $H$-subfields of $K^*$ that is the identity on $K(C^*)^\ell$ and sends $a_1$ to $a_2$.
By Proposition~\ref{prop:ntexpclosure}, we may extend this to a logarithmic $H$-field isomorphism of models of $T^{\nt}_{\exp}$.
Then this is an $\cL_{\elem}$-isomorphism by~\ref{trigmainthmiii} and Lemma~\ref{lem:trigembeddingprep}, so it is elementary by~\ref{trigmainthmi}.

The proof of~\ref{trigmainthmv} follows from the proof of Corollary~\ref{cor:stablyembedded} in much the same way:\ take the definable closures $C_1,C_2$ there in $C^*$ as a model of $\Th(\R_{\elem})$, so that the isomorphism $\iota\colon C_1\to C_2$ also respects the restricted trigonometric functions, pass to the Newton-exponential closure of $K(C_1)^\ell$, and apply~\ref{trigmainthmi},~\ref{trigmainthmiii}, and Lemma~\ref{lem:trigembeddingprep}.
\end{proof}

\section{Transseries, Hardy fields, and surreal numbers}\label{sec:examples}
In this section, we discuss the most important models of $T^{\nt}_{\exp}$.

\subsection{Transseries}\label{subsec:transseries}
Let $\T_{\exp}$ be the natural expansion of the field of logarithmic-exponential transseries to an $\cL_{\log}$-structure.
In the introduction, $\T_{\exp}$ was the expansion of $\T$ by the exponential function and here, it is the expansion by the logarithm function, but there is no difference in terms of model completeness or definability.

The logarithm on $\T_{\exp}$ is a group isomorphism $\T^>\to \T$ that extends the real logarithm and satisfies~\ref{L1} and~\ref{L2}; see~\cite{DMM97} or~\cite[Appendix A]{ADH17} (but note that $\T_{\exp}$ denotes the subfield of exponential transseries in that appendix).
As $\T$ is an $\upo$-free newtonian $H$-field with small derivation by~\cite[Corollary~15.0.2]{ADH17}, we have the following:

\begin{corollary}\label{cor:maincor}
The $\cL_{\log}$-theory of $\T_{\exp}$ is model complete and completely axiomatized by $T^{\nt}_{\exp,\dsmall}$.
\end{corollary}

Many important $H$-subfields of $\T_{\exp}$ are logarithmic $H$-subfields.
To help us find examples, we use the following criterion:

\begin{lemma}\label{lem:Tsubfields}
Let $\log$ denote the logarithm on $\T_{\exp}$ and let $K\supseteq \R$ be an $H$-subfield of $\T_{\exp}$.
If $K$ is henselian and $\log(K^>)\subseteq K$, then $K$ expanded by $\log|_{K^>}$ is a model of $T_{\log}$.
If $K$ is Liouville closed, then $\log(K^>) = K$.
If $K \models T^{\nl}$, then $K$ expanded by $\log|_{K^>}$ is a model of $T^{\nt}_{\exp,\dsmall}$.
\end{lemma}
\begin{proof}
Since the axioms in Definition~\ref{def:logarithm} are universal and $C = \R$, if $K$ is henselian and closed under $\log$, then $K\models T_{\log}$.
Now assume that $K$ is Liouville closed, so in particular $K$ is henselian.
Let $a \in K$ and take $f \in K^>$ with $f^\dagger = a'$, so $(\exp a)/f \in \R^>$ and $\exp a\in K$.
Likewise, for $a\in K^>$, taking $g\in K$ with $g'=a^\dagger$ gives $\log a - g\in\R$, so $\log a \in K$.
\end{proof}

By Lemma~\ref{lem:Tsubfields}, the subfield $\T^{\operatorname{da}}$ of transseries that are $\d$-algebraic over $\Q$ is an $\upo$-free newtonian exponential $H$-subfield of $\T_{\exp}$, as is the subfield $\T_{\operatorname{g}}$ of grid-based transseries.
Thus, $\T^{\operatorname{da}}$ and $\T_{\operatorname{g}}$ are both elementary $\cL_{\log}$-substructures of $\T_{\exp}$.
The field $\T_{\log}$ of logarithmic transseries is an $\upo$-free newtonian logarithmic $H$-subfield of $\T_{\exp}$, but not an exponential $H$-subfield.
See~\cite{vdH06} for more on $\T_{\operatorname{g}}$ and~\cite{Geh17D} for more on $\T_{\log}$.

Using Theorem~\ref{thm:completions} and Corollary~\ref{cor:maincor}, we can transfer some known facts about $\T_{\exp}$ into formal consequences of $T^{\nt}_{\exp}$.

\begin{corollary}\label{cor:modelsT}
Let $K \models T^{\nt}_{\exp}$.
Then the underlying ordered exponential field of $K$ is a model of $\Th(\R_{\exp})$.
\end{corollary}
\begin{proof}
If this corollary holds for $K^{\phi}$ where $\phi \in K^>$, then it holds for $K$, so by compositionally conjugating by a suitable element of $K^>$, we may assume that $K\models T^{\nt}_{\exp,\dsmall}$.
The underlying ordered exponential field of $\T_{\exp}$ is a model of $\Th(\R_{\exp})$ by~\cite[Corollary~2.8]{DMM97}, so this is true for $K$ as well, as $K$ is elementarily equivalent to $\T_{\exp}$.
\end{proof}

Let $T$ be an o-minimal theory extending the theory of real closed ordered fields in a suitable language $\cL$.
Let $M\models T$ and let $\der$ be a derivation on $M$.
Following the terminology of~\cite{FK21}, we say that $\der$ is a \textbf{$T$-derivation on $M$} if for any open $U\subseteq M^n$, any continuously differentiable function $f\colon U\to M$ that is $\cL$-definable without parameters, and any $u= (u_1,\ldots,u_n) \in U$, we have
\[
\der f(u) = \nabla f(u)\cdot (\der u_1,\ldots,\der u_n)
\]
where $\nabla f$ is the gradient of $f$.
An \textbf{$H_T$-field}, as introduced in~\cite{Kap23}, is a model $M \models T$ equipped with a $T$-derivation $\der$ such that $(M,\der)$ is an $H$-field.

\begin{corollary}\label{cor:Tderiv}
Let $T \coloneqq \Th(\R_{\exp})$ and let $K \models T^{\nt}_{\exp}$.
Then $\der$ is a $T$-derivation on $K$ and $K$ is an $H_T$-field.
Moreover, the theory $T^{\nt}_{\exp}$ is the model companion of the theory of $H_T$-fields (viewed as $\cL_{\log}$-structures).
\end{corollary}
\begin{proof}
Being a $T$-derivation is invariant under compositional conjugation by elements of $K^>$, so we may assume that $K\models T^{\nt}_{\exp,\dsmall}$.
Then $K$ is elementarily equivalent to $\T_{\exp}$.
The derivation on $\T_{\exp}$ is a $T$-derivation by~\cite[Example~2.16]{FK21} and being a $T$-derivation is expressible as a collection of $\cL_{\log}$-sentences, so $\der$ is a $T$-derivation as well.
By Corollary~\ref{cor:modelsT}, $K$ is an $H_T$-field.

Every $H_T$-field $E$ is in particular a logarithmic $H$-field with $C_E\models T$ by~\cite[Lemma~2.3]{FK21}, and so can be extended to a model of $T^{\nt}_{\exp}$ by Corollary~\ref{cor:upontexpclosure}.
As this model is itself an $H_T$-field by the first part, and $T^{\nt}_{\exp}$ is model complete, the second part of the corollary follows.
\end{proof}

Let $\T_{\elem}$ be the natural expansion of $\T_{\exp}$ to an $\cL_{\elem}$-structure; see~\cite{DMM97} for details.
One can easily check that the restricted trigonometric functions in this expansion satisfy~\ref{RT1}--\ref{RT3}, so $\T_{\elem} \models T^{\nt}_{\elem}$.
In particular, we have:

\begin{corollary}\label{cor:trigmaincor}
The $\cL_{\elem}$-theory of $\T_{\elem}$ is model complete.
\end{corollary}

We also have analogs of Corollaries~\ref{cor:modelsT} and~\ref{cor:Tderiv}.

\begin{corollary}\label{cor:elemTderiv}
Let $T \coloneqq \Th(\R_{\elem})$ and let $K \models T^{\nt}_{\elem}$.
Then the expansion of the underlying ordered field of $K$ by $\exp$, $\sin$, and $\cos$ is a model of $T$, the derivation $\der$ is a $T$-derivation on $K$, and $K$ is an $H_T$-field.
Moreover, the theory $T^{\nt}_{\elem}$ is the model companion of the theory of $H_T$-fields (viewed as $\cL_{\elem}$-structures).
\end{corollary}
\begin{proof}
As in the proofs of Corollaries~\ref{cor:modelsT} and~\ref{cor:Tderiv}, we may assume that $K$ has small derivation.
Then $K$ is elementarily equivalent to $\T_{\elem}$ by Theorem~\ref{thm:trigmainthm}\ref{trigmainthmii}.
By~\cite[Corollary~2.8]{DMM97} and~\cite[Example~2.16]{FK21}, $\T_{\elem}$ expands a model of $T$ and the derivation on $\T$ is a $T$-derivation, so the same holds for $K$.
In particular, $K$ is an $H_T$-field.

For the second part, let $E$ be an $H_T$-field.
Then $C_E \models T$ by~\cite[Lemma~2.3]{FK21} and $E$ is a logarithmic $H$-field, so Corollary~\ref{cor:upontexpclosure} gives an extension $M\models T^{\nt}_{\exp}$ of $E$ with $C_M = C_E$.
By Theorem~\ref{thm:trigmainthm}\ref{trigmainthmiii}, $M$ admits a unique expansion to a model of $T^{\nt}_{\elem}$ extending the $\cL_{\elem}$-structure on $C_E$ induced by $E$.
Arguing as in the proof of Lemma~\ref{lem:trigembeddingprep}, one can check that this expansion makes $M$ an $\cL_{\elem}$-extension of $E$.
Thus every $H_T$-field extends to a model of $T^{\nt}_{\elem}$, so the second part of the corollary follows from the first part and Theorem~\ref{thm:trigmainthm}\ref{trigmainthmi}.
\end{proof}

\subsection{Hardy fields}\label{subsec:hardy}
Recall that a Hardy field is a set of germs of real-valued functions at $+\infty$ that is closed under differentiation and that forms a field under addition and multiplication.
Every Hardy field containing $\R$ (that is, containing the germs of all constant functions) is an $H$-field with constant field $\R$.

\begin{fact}[{\cite[Chapitre~V]{Bou76}}]\label{fact:liouvillehardy}
Any maximal Hardy field is real closed, contains $\R$, and is closed under taking integrals, exponentials, and logarithms of positive elements.
\end{fact}

\begin{fact}[\cite{ADH24}]\label{fact:uponewthardy}
Any maximal Hardy field is $\upo$-free and newtonian.
\end{fact}

\begin{lemma}\label{lem:loghardy}
Let $H$ be a henselian Hardy field containing $\R$, and suppose $\log f \in H$ for each $f \in H^>$.
Then $H$ equipped with the logarithm $f\mapsto\log f$ is a model of $T_{\log}$.
\end{lemma}
\begin{proof}
\ref{L1} holds by Taylor expansion and~\ref{L2} is just the chain rule.
\end{proof}

Combining Facts~\ref{fact:liouvillehardy} and~\ref{fact:uponewthardy} with Lemma~\ref{lem:loghardy}, we have:

\begin{corollary}\label{cor:maximalhardy}
Any maximal Hardy field is a model of $T^{\nt}_{\exp,\dsmall}$.
In particular, all maximal Hardy fields are elementarily equivalent to $\T_{\exp}$ as $\cL_{\log}$-structures.
\end{corollary}

This can be leveraged to prove the following transfer theorem:

\begin{theorem}\label{thm:transfer}
Let $H$ be a henselian Hardy field containing $\R$ and closed under taking logarithms.
Let $\iota\colon H\to\T_{\exp}$ be a logarithmic $H$-field embedding.
Let $\sigma$ be an $\cL_{\log}(H)$-sentence, and let $\iota(\sigma)$ be the $\cL_{\log}(\T_{\exp})$-sentence obtained by replacing each parameter in $\sigma$ by its image under $\iota$.
Then the following are equivalent:
\begin{enumerate}[label=\normalfont(\arabic*), align = right]
\item $M\models\sigma$ for some maximal Hardy field $M\supseteq H$;
\item $M\models\sigma$ for every maximal Hardy field $M\supseteq H$;
\item $\T_{\exp}\models\iota(\sigma)$.
\end{enumerate}
\end{theorem}
\begin{proof}
Fix a maximal Hardy field $M\supseteq H$, viewed as a model of $T^{\nt}_{\exp,\dsmall}$, so $H$ is a logarithmic $H$-subfield of $M$.
We need to show that $M\models\sigma$ if and only if $\T_{\exp}\models\iota(\sigma)$.
First suppose $H = \R$ and take an $|M|^+$-saturated elementary extension $L$ of $\T_{\exp}$.
By Corollary~\ref{cor:bigembedding}, the inclusion $\R\subseteq L$ extends to a logarithmic $H$-field embedding $M\to L$, which is elementary since $T^{\nt}_{\exp}$ is model complete, so $M$ and $\T_{\exp}$ have the same $\cL_{\log}(\R)$-theory.

Now suppose that $H\neq \R$.
Then $H$ is nontrivially valued so we can take $f\in H^>$ with $f\succ 1$.
For each $n$, the $n$th iterated logarithm $\log_nf$ lies in $H$, and $\iota(\log_nf) = \log_n\iota(f)$.
The sequence $(\log_n\iota(f))$ is coinitial in the set of positive infinite elements of $\T$~\cite[Appendix~A]{ADH17}, so the image of $\Gamma_H^<$ in $\Gamma_{\T}$ under $\iota$ is cofinal in $\Gamma_{\T}^<$.
As $\T$ is $\upo$-free, $\iota(H)$ is as well by Fact~\ref{fact:dalgupofree}, so $H\models T^{\upo}_{\log}$.
Thus, $M$ and $\T_{\exp}$ have the same $\cL_{\log}(H)$-theory by Theorem~\ref{thm:mainthm}.
\end{proof}

\begin{remark}\label{rem:hardytrig}
Let $M$ be a maximal Hardy field.
For $f\in M$ with $f\preceq 1$, the germs $\cos f$ and $\sin f$ lie in $M$ by~\cite[Proposition~6.9]{ADH26}.
Thus $M$ admits restricted trigonometric functions extending the restricted sine and cosine on $\R$ and satisfying~\ref{RT1}--\ref{RT3}.
Expanding $M$ by these functions, we get a model of $T^{\nt}_{\elem}$ with small derivation.
In particular, all maximal Hardy fields are elementarily equivalent to $\T_{\elem}$ as $\cL_{\elem}$-structures, by Theorem~\ref{thm:trigmainthm} and Corollary~\ref{cor:trigmaincor}.
\end{remark}

\subsection{Surreal numbers}\label{subsec:surreal}
The field $\No$ of surreal numbers is a proper class-sized real closed field extension of $\R$ introduced by Conway~\cite{Con76}.
The surreals admit a (surjective) logarithm, defined by Kruskal and Gonshor~\cite{Gon86}, and with this logarithm, $\No$ is an elementary extension of $\R_{\exp}$~\cite{DE01}.
More recently, Berarducci and Mantova equipped the surreals with a derivation that makes $\No$ an exponential $H$-field with small derivation and real exponential constant field $\R$~\cite{BM18}.
With this derivation, the $H$-field $\No$ is newtonian and $\upo$-free~\cite{ADH19}, so $\No$ is a model of $T^{\nt}_{\exp,\dsmall}$.

Each surreal number has an ordinal-valued length, and for each ordinal $\gamma$, we let $\No(\gamma)$ be the set of surreal numbers of length $<\gamma$.
For our purposes, we only need that for $\kappa$ a regular uncountable cardinal, the set $\No(\kappa)$ is an $\upo$-free newtonian Liouville closed $H$-subfield of $\No$ containing $\R$~\cite[Corollary~4.6]{ADH19}.
By adapting the proof of Lemma~\ref{lem:Tsubfields}, we see that $\No(\kappa)$ is closed under $\exp$ and $\log$, so it is an elementary $\cL_{\log}$-substructure of $\No$.

The next proposition, an analog of~\cite[Theorem~3]{ADH19}, shows that the surreal numbers are \emph{universal} among models of $T_{\log}$ with small derivation and archimedean constant field.

\begin{proposition}\label{prop:Noembedding}
Let $K$ be a set-sized model of $T_{\log}$ with archimedean constant field and small derivation.
Then $K$ admits a logarithmic $H$-field embedding into $\No$.
\end{proposition}
\begin{proof}
It suffices to show that some logarithmic $H$-field extension of $K$ with the same constant field admits such an embedding into $\No$.
Since $K$ has small derivation, either $0 \in (\Gamma^<)'$ or $K$ has gap 0.
In the second case, we can use Lemma~\ref{lem:gap2wfree} to extend $K$ to an $\upo$-free logarithmic $H$-field with $0 \in (\Gamma^<)'$, so by passing to an extension, we may assume that $0 \in (\Gamma^<)'$.
Then any logarithmic $H$-field extension of $K$ has small derivation.
By Corollary~\ref{cor:upontexpclosure}, $K$ has an $\upo$-free newtonian exponential $H$-field extension with the same constant field, so we may assume that $K \models T^{\nt}_{\exp,\dsmall}$.
Let $\kappa\coloneqq |K|^+$.
Then $\No(\kappa)\models T^{\nt}_{\exp,\dsmall}$ and by~\cite[Lemma~5.3]{ADH19}, the underlying ordered sets of $\No(\kappa)$ and $\Gamma_{\No(\kappa)}$ are $\kappa$-saturated.
Since $C$ is an archimedean model of $\Th(\R_{\exp})$, the Laskowski--Steinhorn theorem~\cite{LS95} gives an elementary ordered exponential field embedding $C\to \R_{\exp}$.
By Corollary~\ref{cor:bigembedding}, this embedding extends to a logarithmic $H$-field embedding $K\to\No(\kappa)$.
\end{proof}

\begin{corollary}\label{cor:HardyNo}
Let $H$ be a Hardy field with $\log f \in H$ for all $f \in H^>$.
Then there is an embedding $H\to \No$ of ordered valued differential fields that preserves logarithms.
\end{corollary}
\begin{proof}
Let $M$ be a maximal Hardy field extension of $H$.
Then $M \models T^{\nt}_{\exp,\dsmall}$ by Corollary~\ref{cor:maximalhardy}, so Proposition~\ref{prop:Noembedding} gives a logarithmic $H$-field embedding $M\to\No$.
The restriction of this embedding to $H$ is the desired embedding.
\end{proof}

As described in~\cite[Proposition~5.1]{ADH19}, there is a unique strongly $\R$-linear ordered exponential field embedding $\iota\colon \T\to \No$ that sends $x$ to $\omega$.
It was shown in~\cite[Theorem~1]{ADH19} that $\iota$ is an elementary embedding of ordered differential fields, and the question was raised there (and in~\cite{BEK16}) as to whether $\iota$ remains elementary after adding the exponential.
As $\T_{\exp}$ and $\No$ are both models of $T^{\nt}_{\exp}$, Theorem~\ref{thm:mainthm} allows us to answer this question:

\begin{corollary}\label{cor:elemTNo}
The embedding $\iota\colon\T\to \No$ is an elementary embedding of ordered differential exponential fields.
\end{corollary}

It is easy to see that $\iota$ is elementary as an $\cL_{\elem}$-embedding as well.

\section{Future Directions}\label{sec:open}
\subsection*{Quantifier elimination}
Perhaps the most pressing next step is to prove a quantifier elimination result for $T^{\nt}_{\exp}$.
As in the case of $T^{\nl}$, quantifier elimination will almost surely have to involve some additions to our language.
Let $K \models T^{\nt}_{\exp}$ and let $\cL_{\log}^{\Upl, \Upo}$ be the extension of $\cL_{\log}$ by the following symbols:
\begin{enumerate}[label=\normalfont(\arabic*), align = right]
\item\label{langinv} A unary function symbol $\imath$, interpreted in $K$ as multiplicative inversion (with $\imath(0)\coloneqq 0$);
\item\label{langUpl} A unary predicate $\Upl$, interpreted in $K$ as the set $\{-y^{\dagger\dagger}:y \in K,\ y \succ 1\}$;
\item\label{langUpo} A unary predicate $\Upo$, interpreted in $K$ as the set
\[
\{a \in K:4y''+ay = 0\text{ for some }y \in K^{\times}\};
\]
\item\label{langdf} A function symbol $\tilde{f}$ for each function $f\colon \R^n\to \R$ definable without parameters in $\R_{\exp}$, interpreted as a function $C^n\to C$ in the natural way (using that $C$ is a model of $\Th(\R_{\exp})$ and taken to be identically zero off of $C^n$);
\item\label{langst} A unary function symbol $\st$, interpreted in $K$ as the standard part map $K\to C$, so $a-\st(a) \in \smallo$ for $a \in \cO$, and $\st(a)\coloneqq 0$ for $a \notin \cO$.
\end{enumerate}

\begin{question}
Does $T^{\nt}_{\exp}$ eliminate quantifiers in the language $\cL_{\log}^{\Upl, \Upo}$?
\end{question}

The predicates $\Upo$ and $\Upl$ were introduced in~\cite[Chapter~16]{ADH17}, and their addition is necessary to prove quantifier elimination for the theory $T^{\nl}$~\cite[Propositions~16.5.3 and~16.5.4]{ADH17}.
They will almost certainly be required to prove quantifier elimination for $T^{\nt}_{\exp}$.

Since $\Th(\R_{\exp})$ has definable Skolem functions, the function symbols $\tilde{f}$, the standard part map, and multiplicative inversion \emph{almost} ensure that $T_{\log}$ is a universal theory:\ if $L \models T_{\log}$ is expanded by these functions and $E \subseteq L$ is a substructure, then $E$ is an $H$-subfield of $L$ that is closed under $\log$ and $C_E$ is an elementary substructure of $C_L$ as an ordered exponential field.
As $E$ may fail to be henselian, it need not be a model of $T_{\log}$, but its extension $E^\ell$ (identified with an $H$-subfield of $L$) is a model of $T_{\log}$.

It may be possible to replace the function symbols $\tilde{f}$ by any language in which $\Th(\R_{\exp})$ has quantifier elimination.
Note that this theory does \emph{not} eliminate quantifiers in the natural language of ordered exponential fields~\cite{vdD84}.
It may also be possible to remove the standard part map, in exchange for interpreting the function symbols $\tilde{f}$ as functions $K^n\to K$, using that $K \models \Th(\R_{\exp})$ by Corollary~\ref{cor:modelsT}.
This interpretation is more in line with the setting of $H_T$-fields~\cite{Kap23}.
In this setting, substructures of $K$ may no longer be logarithmic $H$-fields, but they can be extended to logarithmic $H$-fields by way of the \emph{logarithmic $H_T$-field hull}~\cite[Proposition~6.16]{Kap23}, where $T$ is the theory of the real field with restricted exponentiation.

A quantifier elimination result would be the first step in analyzing the combinatorial complexity of $T^{\nt}_{\exp}$.
The theory $T^{\nl}$ is distal~\cite[Corollary~7.3]{ACGZ22}; in particular, it has NIP~\cite[Proposition~16.6.6]{ADH17}.
Both proofs use quantifier elimination, and while it seems likely that $T^{\nt}_{\exp}$ is also distal, this is difficult to show without QE.

Quantifier elimination should also allow us to improve Theorem~\ref{thm:transfer} by showing that any two maximal Hardy fields are elementarily equivalent as logarithmic $H$-fields over any common Hardy subfield, without assuming that this subfield embeds into $\T_{\exp}$.
This was shown without the logarithm in~\cite[Theorem~12.3]{ADH24}.

\subsection*{A dimension theory}
A dimension theory for $\T$, and more generally for models of $T^{\nl}$, was developed in~\cite{ADH17B}.
Let $K\models T^{\nl}$ and let $X \subseteq K^n$ be nonempty and definable.
The \textbf{topological dimension} of $X$ is the largest $m$ such that the image of $X$ under some coordinate projection $K^n\to K^m$ has nonempty interior.
It was shown that the differential dimension of $X$ (in the sense of~\cite[Section~1]{ADH17B}) equals its topological dimension~\cite[Corollary~3.2]{ADH17B}, and that $X$ has dimension zero if and only if $X$ is discrete if and only if $X$ is co-analyzable relative to the constant field~\cite[Propositions~4.1 and~6.2]{ADH17B}.

\begin{question}
Is there a good dimension theory for models of $T^{\nt}_{\exp}$? For $K \models T^{\nt}_{\exp}$ and $X\subseteq K^n$ definable and nonempty, are the following equivalent?
\begin{enumerate}[label=\normalfont(\arabic*), align = right]
\item\label{dim1} $X$ has topological dimension zero;
\item\label{dim2} $X$ has differential dimension zero;
\item\label{dim3} $X$ is discrete;
\item\label{dim4} $X$ is co-analyzable relative to the constant field.
\end{enumerate}
\end{question}

Most of these equivalences follow as in~\cite{ADH17B}:\ the equivalence of~\ref{dim1} and~\ref{dim3} can be established using Corollary~\ref{cor:lomin}, as in the proof of~\cite[Proposition~4.1]{ADH17B}, and the equivalence of~\ref{dim2} and~\ref{dim4} follows from the proof of~\cite[Proposition~6.2]{ADH17B}.
We have that~\ref{dim2} implies~\ref{dim1} by~\cite[Lemma~4.4.10]{ADH17}, so the only implication missing is that~\ref{dim1} implies~\ref{dim2}.
As in~\cite{ADH17B}, this requires a better understanding of definable sets, likely via quantifier elimination as discussed above.

\subsection*{Initial embeddings into the surreals}
Let $\iota\colon\T\to\No$ be as in Corollary~\ref{cor:elemTNo}.
In~\cite[Corollary~12.2]{EK21}, $\iota$ was shown to also be \emph{initial} with respect to the binary tree structure on $\No$.
The following question is an exponential analog of a question asked in the final remarks of~\cite{ADH19}.

\begin{question}
Does every model of $T^{\nt}_{\exp,\dsmall}$ with constant field $\R$ admit an initial logarithmic $H$-field embedding into $\No$ over $\R$?
\end{question}

We note that every elementary extension of $\R_{\elem}$ is isomorphic to an initial subfield of $\No$ by~\cite[Corollary~10.5]{EK21}.
Recently, Freni generalized this to elementary extensions of $\R_{\exp}$~\cite{Fre24}.

\subsection*{The exponential DIVP}
Let $K$ be an exponential $H$-field, viewed as an exponential ring with $E\coloneqq \exp$, and let $K[Y_0,\ldots,Y_n]^E$ denote the ring of exponential polynomials over $K$; see~\cite[Section~1]{vdD84B}.
We define the ring of \textbf{differential exponential polynomials} to be
\[
K\{Y\}^E\coloneqq \bigcup_n K[Y_0,\ldots,Y_n]^E,
\]
where we write $Y^{(n)}$ for $Y_n$.
The \textbf{order} of $P \in K\{Y\}^E$ is the least $n$ with $P\in K[Y_0,\ldots,Y_n]^E$.
As in the proof of~\cite[Lemma~3.2]{vdD84B}, there is a unique derivation on $K\{Y\}^E$ that extends the derivation on $K$, sends $Y^{(n)}$ to $Y^{(n+1)}$, and satisfies $\der E(P) = E(P)\der P$.
For $f \in K$, there is a unique differential exponential $K$-algebra morphism $K\{Y\}^E\to K$ that sends $Y$ to $f$ by~\cite[Proposition~1.2]{vdD84B}, and we let $P(f)$ be the image of $P \in K\{Y\}^E$ under this morphism.

As we have already seen, models of $T^{\nt}_{\exp}$ enjoy the DIVP.
We can ask whether this also holds for differential exponential polynomials.
Say that an exponential $H$-field $K$ has $\exp$-DIVP if for all $P\in K\{Y\}^E$ and all $f<g \in K$ with $P(f)<0<P(g)$, there is $y \in K$ with $f<y<g$ and $P(y) = 0$.

\begin{question}
Does every model of $T^{\nt}_{\exp}$ have $\exp$-DIVP?
\end{question}

As models of $T^{\nt}_{\exp}$ expand models of $\Th(\R_{\exp})$ by Corollary~\ref{cor:modelsT}, the underlying ordered exponential field of a model of $T^{\nt}_{\exp}$ is o-minimal, so it has the $\exp$-DIVP for differential exponential polynomials of order 0.
Using Theorem~\ref{thm:mainthm}, Corollary~\ref{cor:Tderiv}, an o-minimal implicit function and cell decomposition argument (see~\cite{AD05,vdD00}), and~\cite[Proposition~8.3 and Remark~8.4]{Kap23}, one should be able to show that every model of $T^{\nt}_{\exp}$ has the $\exp$-DIVP for differential exponential polynomials of order 1 as well.

\bibliographystyle{amsplain}
\bibliography{Texp}

\end{document}